\documentclass[reqno,centertags,12pt]{amsart}
\usepackage{rjmacros4}
\begin{document}
\title[The peak model for supersingular perturbations]{The peak model
for finite rank supersingular perturbations}
\author{Rytis Jur\v{s}\.{e}nas}
\address{Vilnius University,
Institute of Theoretical Physics and Astronomy,
Saul\.{e}tekio ave.~3, 10257 Vilnius, Lithuania}
\email{rytis.jursenas@tfai.vu.lt}
\thanks{The research was inspired by the
topics presented at Insubria Summer School in
Mathematical Physics, University of Insubria, Como (Italy),
18-22 September 2017.}
\date{\today}
\begin{abstract}
In its original form
the peak model for rank one supersingular perturbations
of class $\fH_{-4}$ or higher of a nonnegative self-adjoint
operator requires that the Gram matrix of the model
should be diagonal.
Here we remove the restriction on the Gram matrix.
In particular we explain the origin of the Krein
$Q$-function associated with the Gram matrix.
\end{abstract}
\maketitle
%%%%%%%%%%%%%%%%%%%%%%%%%%%%%%%%%%%%%%%%%%%%%%%%%%%%%%%%%%%%%%
%%%%%%%%%%%%%%%%%%%%%%%%%%%%%%%%%%%%%%%%%%%%%%%%%%%%%%%%%%%%%%
%%%%%%%%%%%%%%%%%%%%%%%%%%%%%%%%%%%%%%%%%%%%%%%%%%%%%%%%%%%%%%
%%%%%%%%%%%%%%%%%%%%%%%%%%%%%%%%%%%%%%%%%%%%%%%%%%%%%%%%%%%%%%
\section{Introduction}
%%%%%%%%%%%%%%%%%%%%%%%%%%%%%%%%%%%%%%%%%%%%%%%%%%%%%%%%%%%%%%
%%%%%%%%%%%%%%%%%%%%%%%%%%%%%%%%%%%%%%%%%%%%%%%%%%%%%%%%%%%%%%
%%%%%%%%%%%%%%%%%%%%%%%%%%%%%%%%%%%%%%%%%%%%%%%%%%%%%%%%%%%%%%
%%%%%%%%%%%%%%%%%%%%%%%%%%%%%%%%%%%%%%%%%%%%%%%%%%%%%%%%%%%%%%
The theory of higher order singular or else supersingular
perturbations of a
self-adjoint operator in a Hilbert space is,
in principal, the theory of generalized, that is,
exit space self-adjoint extensions,
where perturbations are interpreted
by means of generalized resolvents or, equivalently,
generalized Nevanlinna families.
For finite rank perturbations, the exit space
$\cH=\fH_m\dsum\fK$ is
made by extending a reference Hilbert space
$\fH_m$ by a disjoint finite-dimensional linear space
$\fK$.
Depending on the precise definition of $\fK$,
the cascade (A and B) and the peak models
for supersingular perturbations are considered
among researchers. Specifically,
the cascade models for rank one
supersingular perturbations of
a nonnegative self-adjoint operator are
developed in \cite{Dijksma05}, see also references therein.
In the B-model $\cH$ is a Pontryagin space
with a nontrivial index of indefiniteness.
In the A-model, whether or not $\cH$ is a Hilbert
space depends on how one defines
the scalar product in the scale of Hilbert spaces
\cite[Theorem~3.2]{Dijksma05}; \cf \cite{Jursenas21}.
For classical, that is, singular
perturbations the reader may consult the monographs
\cites{Simon05a,Albeverio00} and references there.

The peak model, first introduced in \cite{Kurasov09},
is our main object of interest here.
In the present model $\cH$ is a Hilbert space,
because the singular elements that span $\fK$
form the Gram matrix, $\cG$, which is positive;
this is the main motive in \cite{Kurasov09}
for introducing an alternative to the cascade models.
On the other hand,
the peak model, as stated, has limitations
in that the elements that generate $\fK$
must be orthogonal or, equivalently, $\cG$
must be diagonal.
A quick example of a $6$-dimensional Laplace operator
with Dirac distribution
already shows that the orthogonality condition
does not necessarily hold.

We recall that, for a $\nu$-dimensional
$(\nu\geq4)$ Laplace operator with Dirac distribution
$m=(\nu-2)/2$ if $\nu$ is even
and $m=(\nu-3)/2$ if $\nu$ is odd.
Two modifications are, for example, as follows.
By lowering $\nu$ while taking instead the distributional
derivative of Dirac distribution,
the perturbation with $\nu=3$
and $m=1$ is considered in \cite[Sec.~10]{Kurasov09}.
A two-particle Rashba spin-orbit coupled operator with
point-interaction considered in
\cite[Example~3.4]{Jursenas19},
\cite[Sec.~5]{Jursenas18a} shows $m=2$ for $\nu=6$.

In the present study we
remove the restriction on $\cG$ by
considering operator extensions to a subspace
of $\cH$ rather than to the whole $\cH$ or,
what is the same, by considering linear relations in
$\cH$ rather than operator extensions;
the operator part of such a linear relation is
precisely the operator extension to $\cH_1$.
The exit space $\cH_1\dfn\fH_m\dsum\fK_1$ is
therefore constructed by taking
a suitable subspace $\fK_1\subseteq\fK$.
In this way we still work within
the framework of the peak model, because
there is the smallest nontrivial subspace,
$\fKmin$, which is disjoint from $\fH_m$
and such that $\fKmin\subseteq\fK_1\subseteq\fK$.
The existence of $\fKmin$,
but in the cascade A-model,
was first observed in \cite{Jursenas21}.

The characterization of the boundary value space
for the operator extensions to $\cH_1$
should be considered as our first main result
out of the three. The second main result
is that a generalized resolvent
corresponding to the operator extension
in $\cH_1$ parametrized by a self-adjoint
linear relation $\Theta$
is in bijective correspondence with a Nevanlinna
family $z\mapsto(\cC^{-1})^*(\mQ(z)-\Theta)\cC^{-1}$.
The transfer matrix $\cC$ serves for a ``scale parameter''
depending on how one defines the scalar product
in the scale of Hilbert spaces; this is because
we avoid attaching ourselves to a specific
(typically polynomial) definition of scale spaces.
The matrix function $\mQ$ is a Nevanlinna function.
For $\cG$ diagonal, this $\mQ$ coincides with that
in \cite{Kurasov09}, where it is termed
the $Q$-function associated with the Gram
matrix $\cG$. We explain the origin of $\mQ$
by demonstrating that it is
the Weyl function corresponding to a boundary
triple of a certain symmetric operator
in $\fK_1$. In the terminology of \cite{Dijksma04a} we
present a realization for $\mQ$; it is
minimal iff $\fK_1=\fKmin$.
The latter is our third main result.

After introductory Sections~\ref{sec:basic},
\ref{sec:pre},
operator extensions to $\cH_1$ are studied
in Section~\ref{sec:gen}. Specifically,
the characterization of the boundary value space of
extensions, including a realization for $\mQ$,
is presented in Theorem~\ref{thm:mBmax}, and
a generalized resolvent corresponding to a self-adjoint
operator extension to $\cH_1$ is presented in
Theorem~\ref{thm:gen}.
In Section~\ref{sec:fin} we study some properties
of $\mQ$ depending on $\fK_1$.
There is also Appendix~\ref{sec:ibt},
where we use the peak model as
an example for an isometric boundary triple
\cite[Definition~1.9]{Derkach17}.

Throughout,
we present our results using the language of
linear relations in Hilbert spaces
and the theory abstract boundary value spaces
\cites{Behrndt18,Derkach17,Derkach12,Bruning08,Behrndt07,
Derkach06,Derkach91,Azizov89}.
Linear relations are referred to as relations
and operators are identified with their graphs
(single-valued relations).
%%%%%%%%%%%%%%%%%%%%%%%%%%%%%%%%%%%%%%%%%%%%%%%%%%%%%%%%%%%%%%
%%%%%%%%%%%%%%%%%%%%%%%%%%%%%%%%%%%%%%%%%%%%%%%%%%%%%%%%%%%%%%
%%%%%%%%%%%%%%%%%%%%%%%%%%%%%%%%%%%%%%%%%%%%%%%%%%%%%%%%%%%%%%
%%%%%%%%%%%%%%%%%%%%%%%%%%%%%%%%%%%%%%%%%%%%%%%%%%%%%%%%%%%%%%
\section{Some background}\label{sec:basic}
%%%%%%%%%%%%%%%%%%%%%%%%%%%%%%%%%%%%%%%%%%%%%%%%%%%%%%%%%%%%%%
%%%%%%%%%%%%%%%%%%%%%%%%%%%%%%%%%%%%%%%%%%%%%%%%%%%%%%%%%%%%%%
%%%%%%%%%%%%%%%%%%%%%%%%%%%%%%%%%%%%%%%%%%%%%%%%%%%%%%%%%%%%%%
%%%%%%%%%%%%%%%%%%%%%%%%%%%%%%%%%%%%%%%%%%%%%%%%%%%%%%%%%%%%%%
Consider a closed symmetric relation $T$ in a Hilbert
space $\fH$, with equal and finite defect numbers
$(d,d)$. Then (\eg \cites{Behrndt11,Derkach06})
$T$ has an (ordinary) boundary triple
$\Pi_\Gamma=(\bbC^d,\Gamma_0,\Gamma_1)$, where the boundary
operator $\Gamma\dfn(\Gamma_0,\Gamma_1)$
from $T^*$ (the adjoint in $\fH$ of $T$) to $\bbC^{2d}$
is surjective, and moreover the Green identity
holds
(\eg \cite[Definition~3.1]{Derkach17},
\cite[Definition~7.11]{Derkach12})
\[
\braket{f,g^\prime}_\fH-
\braket{f^\prime,g}_\fH=
\braket{\Gamma_0\whf,\Gamma_1\whg}_{\bbC^d}-
\braket{\Gamma_1\whf,\Gamma_0\whg}_{\bbC^d}
\]
for all $\whf=(f,f^\prime)$ and
$\whg=(g,g^\prime)$ from $T^*$. The scalar
product in $\fH$ is denoted by $\braket{\cdot,\cdot}_\fH$.
Equivalently, by considering
$\Gamma$ with domain $T^*$ as an operator
from a Krein space $\fH^2$ with fundamental
symmetry ($I_\fH$ is the identity in $\fH$)
\[
\whJ_\fH\dfn
\begin{pmatrix}
0 & -\img I_\fH \\ \img I_\fH & 0
\end{pmatrix}\co
\begin{matrix}
\fH \\ \op \\ \fH
\end{matrix}\lto
\begin{matrix}
\fH \\ \op \\ \fH
\end{matrix}
\]
and an indefinite metric
\[
[\whf,\whg]\dfn-\img
\braket{f,g^\prime}_\fH+\img\braket{f^\prime,g}_\fH
\]
to a Krein space $\bbC^{2d}$, with the fundamental
symmetry $\whJ_{\bbC^d}$ and an indefinite
metric defined similarly,
$\Pi_\Gamma$ is said to be a boundary triple for
$T^*$ if $\Gamma$ is unitary, \ie if the inverse
$\Gamma^{-1}$ in the sense of relations
coincides with the Krein
space adjoint $\Gamma^+\dfn\whJ_\fH\Gamma^*\whJ_{\bbC^d}$;
$\Gamma^*$ is a Hilbert space adjoint.
See \cite{Azizov89} for Krein spaces.

Associated with $\Pi_\Gamma$ is the
$\gamma$-field $\gamma_\Gamma$ and the matrix
valued Weyl function $M_\Gamma$ defined by
\[
\begin{split}
&\gamma_\Gamma(z)\dfn\pi_1\whgm_\Gamma(z)\,,\quad
M_\Gamma(z)\dfn\Gamma_1\whgm_\Gamma(z)\,,
\quad z\in\bbC\setm\bbR\,,
\\
&\whgm_\Gamma\dfn(\Gamma_0\vrt_{\whfN_z(T^*)})^{-1}\,,
\quad
\pi_1\co\fH^2\lto\fH\,,\quad
\whf\mapsto f\,.
\end{split}
\]
The defect subspaces
\[
\fN_z(T^*)\dfn\ker(T^*-z)\,,\quad
\whfN_z(T^*)\dfn\{(f,zf)\vrt
f\in\fN_z(T^*) \}\,.
\]
We use the notation $\dom T$, $\ker T$,
$\ran T$ for the domain, kernel, range of $T$.
As a rule we omit the identity operator.
Because $T_0\dfn\ker\Gamma_0$ is a self-adjoint
relation in $\fH$, by the von Neumann formula
the functions $\gamma_\Gamma$
and $M_\Gamma$ extend to $\rho(T_0)$,
the resolvent set of $T_0$.
If $T$ is densely defined, $\Gamma$ on
$T^*$ is identified with $\Gamma$ on
$\dom T^*$, as well as $\whgm_\Gamma$ is
identified with $\gamma_\Gamma$.

The Weyl function $M_\Gamma$ corresponding to
a boundary triple $\Pi_\Gamma$ for $T^*$ is both
a Krein $Q$-function for a pair
$(T,T_0)$
\cites{Arlinski15,Derkach99}, \ie it satisfies
\begin{equation}
M_\Gamma(z)-M_\Gamma(z_0)^*=(z-\ol{z_0})
\gamma_\Gamma(z_0)^*\gamma_\Gamma(z)\,,
\quad z,z_0\in\rho(T_0)
\label{eq:KQ}
\end{equation}
and a Nevanlinna function
\cite{Dijksma18}, \ie it is analytic on
$\bbC\setm\bbR$ and satisfies
\[
M_\Gamma(z)^*=M_\Gamma(\ol{z})\,,\quad
\Im M_\Gamma(z)/\Im z\geq0\,.
\]
Nevanlinna function is a special case
of a Nevanlinna family, defined analogously
but for relations. See \cite{Derkach06}
for various subclasses.

Let $\wtT$ be a self-adjoint extension
of $T$ in some possibly larger Hilbert
space $\wtfH\supseteq\fH$. Let
$P_\fH$ be an orthogonal projection
in $\wtfH$ onto $\fH$.
There is
(\eg \cite[Theorems~6.1, 6.2]{Derkach09})
a unique $d\times d$ relation valued
Nevanlinna family $\tau$ such that
\begin{align*}
P_\fH(\wtT-z)^{-1}\vrt_\fH=&(T_{-\tau(z)}-z)^{-1}
\\
=&(T_0-z)^{-1}-
\gamma_\Gamma(z)(\tau(z)+M_\Gamma(z))^{-1}
\gamma_\Gamma(\ol{z})^*
\end{align*}
for $z\in\rho(\wtT)\cap\rho(T_0)$.
The above generalized Krein--Naimark resolvent formula
establishes a bijective correspondence between
the sets of all generalized resolvents
$P_\fH(\wtT-z)^{-1}\vrt_\fH$ of $T$
and all $d\times d$ relation valued Nevanlinna families
$\tau$;
\cites{Dijksma18,Derkach09,Langer77}.
Particularly, $z\mapsto\tau(z)$ is a matrix function
iff $\wtT\cap T_0=T$,
while $\tau(z)\equiv-\Theta$ is constant iff
$\wtT\in\Ext(T)$, \ie $T\subseteq\wtT\subseteq T^*$.
The \v{S}traus family $z\mapsto T_{-\tau(z)}$ in
$\fH$ corresponding to $\wtT$ is
given by $z\mapsto\ker(\Gamma_1+\tau(z)\Gamma_0)$,
see also \cite[Theorem~2.7.3]{Behrndt20}.
%%%%%%%%%%%%%%%%%%%%%%%%%%%%%%%%%%%%%%%%%%%%%%%%%%%%%%%%%%%%%%
%%%%%%%%%%%%%%%%%%%%%%%%%%%%%%%%%%%%%%%%%%%%%%%%%%%%%%%%%%%%%%
%%%%%%%%%%%%%%%%%%%%%%%%%%%%%%%%%%%%%%%%%%%%%%%%%%%%%%%%%%%%%%
%%%%%%%%%%%%%%%%%%%%%%%%%%%%%%%%%%%%%%%%%%%%%%%%%%%%%%%%%%%%%%
\section{Preparatory statements}\label{sec:pre}
%%%%%%%%%%%%%%%%%%%%%%%%%%%%%%%%%%%%%%%%%%%%%%%%%%%%%%%%%%%%%%
%%%%%%%%%%%%%%%%%%%%%%%%%%%%%%%%%%%%%%%%%%%%%%%%%%%%%%%%%%%%%%
%%%%%%%%%%%%%%%%%%%%%%%%%%%%%%%%%%%%%%%%%%%%%%%%%%%%%%%%%%%%%%
%%%%%%%%%%%%%%%%%%%%%%%%%%%%%%%%%%%%%%%%%%%%%%%%%%%%%%%%%%%%%%
Throughout, $L$ is a self-adjoint operator
in a (complex and separable) Hilbert space $\fH_0$
with scalar product $\braket{\cdot,\cdot}_0$
and norm $\norm{\cdot}_0\dfn\sqrt{\braket{\cdot,\cdot}_0}$.
It is not assumed that $L$ is semibounded.
\begin{defn}
A sequence of Hilbert spaces
\[
\cdots\subseteq\fH_2\subseteq\fH_1\subseteq
\fH_0\subseteq\fH_{-1}\subseteq\fH_{-2}\subseteq
\cdots
\]
is said to be the scale of Hilbert spaces
associated with $L$, and each $\fH_n$ taken separately
is termed the scale space, if the following
conditions hold for $n\geq0$:
\begin{SL}
\item[(a)]
$\fH_n=(\dom L^{n/2},\braket{\cdot,\cdot}_n)$
with scalar product
\[
\braket{f,g}_n\dfn
\braket{\Omega_nf,\Omega_ng}_0\,,\quad
f,g\in\fH_n\,,
\]
the induced norm
$\norm{f}_n\dfn\sqrt{\braket{f,f}_n}$,
and a unitary operator
$\Omega_n$ from $\fH_n$ onto $\fH_0$,
$\Omega_0\dfn I$ (identity), such that:
\begin{SL}
\item[(a$_1$)]
$\Omega_1(\fH_2)=\fH_1$ and
$\Omega_n(\fH_{n+2})=\fH_2$,
\item[(a$_2$)]
$L\Omega_n=\Omega_nL$,
\item[(a$_3$)]
$\Omega_n$ is self-adjoint in $\fH_0$,
\item[(a$_4$)]
the $\Omega_n$'s
are mutually commuting.
\end{SL}
\item[(b)]
The strong dual $\fH_{-n}$ of $\fH_n$
is a Hilbert space with the scalar product,
$\braket{\cdot,\cdot}_{-n}$, defined
via the duality pairing $\braket{\cdot,\cdot}\co
\fH_{-n}\times\fH_n\lto\bbC$ as follows:
\[
\braket{\psi,\phi}_{-n}\dfn
\braket{\wtilde{\Omega}^{-1}_n\psi,
\wtilde{\Omega}^{-1}_n\phi}_0\,,\quad
\psi,\phi\in\fH_{-n}
\]
where a unitary operator $\wtilde{\Omega}_n$ from
$\fH_0$ onto $\fH_{-n}$ is defined by
\[
\braket{\wtilde{\Omega}_nu,f}\dfn
\braket{u,\Omega_nf}_0\,,\quad
u\in\fH_0\,,\quad f\in\fH_n\,.
\]
\end{SL}
\end{defn}
For $n\geq0$,
$L_n$ is the domain restriction to
$\fH_{n+2}$ of $L$. Notice that
the range of $L_n$ is contained in $\fH_n$.
An operator
$L_n$ is self-adjoint in $\fH_n$
iff $\fH_n$ is dense in $\norm{\cdot}_0$.
By (a$_{1,2}$) and the Riesz representation
theorem
\[
\braket{\psi,f}=
\braket{\wtilde{\Omega}^{-1}_n\psi,
\Omega_nf}_0
\]
$L_n$ is self-adjoint in $\fH_n$, so
$\fH_{n+1}\subseteq\fH_n$ densely and continuously.
Similarly
$\fH_{-n}$ is dense in $\fH_{-n-1}$ by (b).
Moreover, the self-adjointness of $L_n$ in $\fH_n$
further yields
\[
\Omega_n(\fH_{n+t})=\fH_{t}
\]
for all nonnegative integers $t$.

The definition of the triplet adjoint,
as stated below, will suffice for our study.
\begin{defn}[\cf \cites{Kurasov09,Dijksma05}]\label{defn:tr}
Consider the triple $\fH_n\subseteq\fH_0\subseteq\fH_{-n}$
from the scale of Hilbert spaces associated
with $L$, and let $L^\prime\subseteq L$ be
an operator in $\fH_n$.
The triplet adjoint of $L^\prime$
with respect to this triple
is a closed relation in $\fH_{-n}$
defined by
\[
\{(\psi,\phi)\in\fH^2_{-n}\vrt
(\forall f\in\dom L^\prime)\;
\braket{\phi,f}=\braket{\psi,Lf} \}\,.
\]
In particular, the triplet adjoint
of $L_n$ is denoted by $L_{-n}$.
\end{defn}
By considering $\Omega_n$ as an operator in
$\fH_0$,
(a$_3$) allows us to view $\wtilde{\Omega}_n\vrt_{\fH_n}$
as $\Omega_n$. Thus, by (a$_{1-4}$)
a bounded operator $L_n$
from a Hilbert space $\fH_{n+2}$ to
a Hilbert space $\fH_n$ has a continuation
$L_{-n}$, which is a self-adjoint operator in $\fH_{-n}$
with dense domain $\fH_{2-n}=\wtilde{\Omega}_n(\fH_2)$.
\begin{rem}
A polynomial description of
$\Omega_n=\sqrt{P_n(L)}$
falls within our definition of the scale space;
$P_n$ is a positive polynomial in $L$ or $\abs{L}$,
of degree $n\geq0$.
Thus:
\begin{SL}
\item[]
$P_n(L)=(\abs{L}+I)^{n}$ in
\cites{Albeverio13,Albeverio10,Albeverio07,Albeverio97,Kurasov03}.
\item[]
$P_n(L)=L^{n}$, $0\in\rho(L)$
in \cite{Derkach03}.
\item[]
$P_n(L)=\prod_{j=1}^n(L-z_j)$,
$L\geq0$, and $z_1,\ldots,z_n<0$
in \cites{Kurasov09,Dijksma05}.
\end{SL}
In these examples
$\wtilde{\Omega}_n=\sqrt{P_n(L_{-n})}$.
\end{rem}
In view of Definition~\ref{defn:tr} and
the preceding remarks, everywhere else below
we omit the index in
$L_n\equiv L$ if no confusion can arise.
With the same meaning $\wtilde{\Omega}_n\equiv\Omega_n$.
%%%%%%%%%%%%%%%%%%%%%%%%%%%%%%%%%%%%%%%%%%%%%%%%%%%%%%%%%%%%%%
%%%%%%%%%%%%%%%%%%%%%%%%%%%%%%%%%%%%%%%%%%%%%%%%%%%%%%%%%%%%%%
%%%%%%%%%%%%%%%%%%%%%%%%%%%%%%%%%%%%%%%%%%%%%%%%%%%%%%%%%%%%%%
%%%%%%%%%%%%%%%%%%%%%%%%%%%%%%%%%%%%%%%%%%%%%%%%%%%%%%%%%%%%%%
\subsubsection*{Symmetric operator}
%%%%%%%%%%%%%%%%%%%%%%%%%%%%%%%%%%%%%%%%%%%%%%%%%%%%%%%%%%%%%%
%%%%%%%%%%%%%%%%%%%%%%%%%%%%%%%%%%%%%%%%%%%%%%%%%%%%%%%%%%%%%%
%%%%%%%%%%%%%%%%%%%%%%%%%%%%%%%%%%%%%%%%%%%%%%%%%%%%%%%%%%%%%%
%%%%%%%%%%%%%%%%%%%%%%%%%%%%%%%%%%%%%%%%%%%%%%%%%%%%%%%%%%%%%%
Fix an integer $m\geq1$ and
consider the family
$\{\vp_s\}^d_{s=1}$ of linearly independent
functionals from $\fH_{-m-2}\setm\fH_{-m-1}$.
In the theory of supersingular perturbations
one looks for a generalized resolvent of
the symmetric restriction, $\Lmin\subseteq L$, to
$f\in\fH_{m+2}$ such that
$\braket{\vp,f}=0$; $\braket{\vp,\cdot}$
stands for the vector valued functional
with components $\braket{\vp_s,\cdot}$.
Then $\Lmin$ is a
closed densely defined symmetric operator in $\fH_m$,
with defect numbers $(d,d)$,
$d<\infty$, and defect subspaces
\[
\begin{split}
&\fN_z(\Lmin^*)=G_z(\bbC^d)\,,\quad
G_z(c)\dfn
\sum_{s}c_sG_s(z)\,,
\quad z\in\rho(L)\,,
\\
&G_s(z)=P(L)^{-1}g_s(z)\,,\quad
g_s(z)\dfn(L-z)^{-1}\vp_s\,,\quad
P(L)\dfn\Omega^2_m
\end{split}
\]
and $c=(c_s)\in\bbC^d$.
Unless specified otherwise, the summation
indexes $s$, $s^\prime$, $\ldots$ run over
$\{1,\ldots,d\}$.

The triplet adjoint $\Lmax$ of $\Lmin$
with respect to $\fH_m\subseteq\fH_0\subseteq\fH_{-m}$
extends $L$ to
$\fH_{2-m}\dsum\fN_z(\Lmax)$ (direct sum,
\cf \cite[Definition~3.1]{Kurasov09})
\[
\fN_z(\Lmax)=g_z(\bbC^d)\,,\quad
g_z(c)\dfn
\sum_{s}c_sg_s(z)\,.
\]
Up to this point one may also look at
\cite{Dijksma05} with $d=1$.
%%%%%%%%%%%%%%%%%%%%%%%%%%%%%%%%%%%%%%%%%%%%%%%%%%%%%%%%%%%%%%
%%%%%%%%%%%%%%%%%%%%%%%%%%%%%%%%%%%%%%%%%%%%%%%%%%%%%%%%%%%%%%
%%%%%%%%%%%%%%%%%%%%%%%%%%%%%%%%%%%%%%%%%%%%%%%%%%%%%%%%%%%%%%
%%%%%%%%%%%%%%%%%%%%%%%%%%%%%%%%%%%%%%%%%%%%%%%%%%%%%%%%%%%%%%
\subsubsection*{Exit space}
%%%%%%%%%%%%%%%%%%%%%%%%%%%%%%%%%%%%%%%%%%%%%%%%%%%%%%%%%%%%%%
%%%%%%%%%%%%%%%%%%%%%%%%%%%%%%%%%%%%%%%%%%%%%%%%%%%%%%%%%%%%%%
%%%%%%%%%%%%%%%%%%%%%%%%%%%%%%%%%%%%%%%%%%%%%%%%%%%%%%%%%%%%%%
%%%%%%%%%%%%%%%%%%%%%%%%%%%%%%%%%%%%%%%%%%%%%%%%%%%%%%%%%%%%%%
Starting from now on, the cascade and the peak models
break apart. In the peak model one considers the restriction,
$\Amax$, of $\Lmax$ to the model space $\cH$,
$\fH_m\subseteq\cH\subseteq\fH_{-m}$, which is
a Hilbert space defined by
\[
\cH\dfn(\fH_m\dsum\fK,\braket{\cdot,\cdot}_\cH)
\]
with an $md$-dimensional Hilbert subspace
$\fK=(\fK,\braket{\cdot,\cdot}_{-m})$
of $\fH_{-m}$ defined by ($\vee$ $\equiv$ linear span)
\[
\fK\dfn\vee\{g_i\dfn g_s(z_j)\}\,,\quad
i=\alpha(s,j)\dfn m(s-1)+j\,.
\]
Clearly the numbers
\[
\cZ\dfn\{z_j\in\rho(L)\}
\]
are assumed distinct.
Unless specified otherwise, the summation
indexes $j$, $j^\prime$, $\ldots$ run over
$\{1,\ldots,m\}$.
\begin{rem}
An index
$i=\alpha(s,j)$ is uniquely determined by $s$
and $j$, that is, the Kronecker symbol
$\delta_{i,i^\prime}=\delta_{s,s^\prime}
\delta_{j,j^\prime}$ for some
$i^\prime=\alpha(s^\prime,j^\prime)$.
\end{rem}
The scalar product
\[
\braket{f+k,f^\prime+k^\prime }_\cH\dfn
\braket{f,f^\prime}_m+
\braket{k,k^\prime}_{-m}
\]
for $f,f^\prime\in\fH_m$ and
$k,k^\prime\in\fK$.
A bijective correspondence
$\fK\ni k\leftrightarrow d(k)\in\bbC^{md}$
is established via the Gram matrix $\cG$
as follows:
\[
\begin{split}
&k=\sum_i d_i(k)g_i\,,\quad
d(k)=(d_i(k))=\cG^{-1}
(\braket{g_i,k}_{-m})\,,
\\
&\cG=(\cG_{ii^\prime}\dfn
\braket{g_i,g_{i^\prime} }_{-m} )\,.
\end{split}
\]
Thus
\[
\braket{k,k^\prime}_{-m}=
\braket{d(k),\cG d(k^\prime) }_{\bbC^{md}}
\]
and in this way
$\cH$ is isomorphic to a Hilbert sum
$\fH_m\op(\bbC^{md},\braket{\cdot,\cG\cdot}_{\bbC^{md}})$.

The set $\fK$ contains a subset,
$\fKmin$,
which is least possible in order the exit space
extensions should cover the case of defect numbers $(d,d)$.
To see this, consider a polynomial
\[
\wtP(L)\dfn\prod_{j}(L-z_j)\,.
\]
It is a bijective operator from
$\fH_{m}$ to $\fH_{-m}$
(with a possible continuation as described
previously), which is seen by using
the sum formula for the inverse
(\cf \cite[Eq.~(6.4)]{Kurasov09}, \cite{Kurasov03b})
\begin{equation}
\wtP(L)^{-1}=\sum_{j}b_j(L-z_j)^{-1}\,,\quad
b_j\dfn\prod_{j\neq j^\prime}
(z_j-z_{j^\prime})^{-1}\,.
\label{eq:sum}
\end{equation}
Thus we have
\begin{lem}
$\fK\supseteq\fKmin$ where
\[
\fKmin\dfn\fK\cap\fH_{m-2}=\vee
\{\ff_s\dfn\wtP(L)^{-1}\vp_s \}\,.
\]
\end{lem}
\begin{proof}
$\vee\{\ff_s\}\subseteq\fKmin$
is due to \eqref{eq:sum}.
Consider an arbitrary $k\in\fKmin$. Since
$\wtP(L)^{-1}(\fH_{-m-2})=\fH_{m-2}$,
$(\exists\psi\in\fH_{-m-2})$
$k=\wtP(L)^{-1}\psi$, and then by \eqref{eq:sum}
\[
0=\sum_j(L-z_j)^{-1}\phi_j\,,\quad
\phi_j\dfn
b_j\psi-\sum_s d_{i}(k)\vp_s\,.
\]
Then
\begin{align*}
0=&\sum_j(L-z_1)(L-z_j)^{-1}\phi_j
\\
=&\sum_{j\geq1}\phi_j+\sum_{j\geq2}
(z_j-z_1)(L-z_2)^{-1}\phi_j
\\
&+\sum_{j\geq3}
(z_j-z_1)(z_j-z_2)(L-z_2)^{-1}
(L-z_3)^{-1}\phi_j
\\
&+\cdots+
(z_m-z_1)\cdots(z_m-z_{m-1})
(L-z_2)^{-1}\cdots(L-z_m)^{-1}\phi_m\,.
\end{align*}
The elements $\sum_{j\geq r}$
belong to mutually disjoint sets
for different $r\in\{1,\ldots,m\}$, namely,
$\fH_{2r-m-4}\setm\fH_{2r-m-3}$, which
shows $(\forall r)\sum_{j\geq r}=0$, and then
$(\forall j)$ $\phi_j=0$. As a result $(\forall j)$
$d_{i}(k)/b_j\equiv \chi_s$ and
$\psi=\sum_s \chi_s\vp_s$, \ie
$k\in\vee\{\ff_s\}$.
\end{proof}
The Gram matrix of $\fKmin$ is denoted by
\[
\cGmin\dfn(\braket{\ff_s,\ff_{s^\prime} }_{-m})
\]
and an element $k=\kmin(\chi)\in\fKmin$ by
\[
\kmin(\chi)\dfn\sum_s\chi_s\ff_s\,,\quad
\chi=(\chi_s)\in\bbC^d\,.
\]
In this way $\fKmin\leftrightarrow\bbC^d$ bijectively
and moreover
\begin{equation}
d(k)=\hb\chi\quad\text{for}\quad
k=\kmin(\chi)\,;\quad
\hb\dfn\cG^{-1}\cG_b\,.
\label{eq:hb}
\end{equation}
The rectangular $md\times d$ matrix
\[
\cG_b\dfn(\braket{g_i,\ff_{s^\prime}}_{-m})\,.
\]
Notice moreover that the matrices $\cGmin$ and $\cG_b$
(or $\hb$) are related by the equality
(with $\cG^*_b$ the adjoint of $\cG_b$)
\[
\cGmin=\cG^*_b\cG^{-1}\cG_b=\cG^*_b\hb\,.
\]

The space $\fKmin$ can be
equivalently represented by using an initially
given $P(L)$ instead of $\wtP(L)$ as follows.
Let
\[
p(L)\dfn P(L)\wtP(L)^{-1}\,,\quad
\fe_s\dfn P(L)^{-1}\vp_s
\]
so that
\[
\begin{split}
&p(L)\fe_s=\ff_s=\sum_{s^\prime}
\cC_{s^\prime s}\fe_{s^\prime}\,,
\quad
\cC=(\cC_{ss^\prime})\dfn\cB^{-1}\cA\,,
\\
&\cA\dfn(\braket{\fe_s,\ff_{s^\prime} }_{-m})\,,
\quad
\cB\dfn(\braket{\fe_s,\fe_{s^\prime} }_{-m})
\;(\text{Gram})\,.
\end{split}
\]
From $\cGmin=\cC^*\cB\cC$ one sees that
$\cC$ (and its adjoint $\cC^*$)
is nonsingular, and
$\fe_s=\sum_{s^\prime}(\cC^{-1})_{s^\prime s}
\ff_{s^\prime}$.
Moreover
\begin{lem}
$\fKmin=\vee\{\fe_s\}$.
\end{lem}
%%%%%%%%%%%%%%%%%%%%%%%%%%%%%%%%%%%%%%%%%%%%%%%%%%%%%%%%%%%%%%
%%%%%%%%%%%%%%%%%%%%%%%%%%%%%%%%%%%%%%%%%%%%%%%%%%%%%%%%%%%%%%
%%%%%%%%%%%%%%%%%%%%%%%%%%%%%%%%%%%%%%%%%%%%%%%%%%%%%%%%%%%%%%
%%%%%%%%%%%%%%%%%%%%%%%%%%%%%%%%%%%%%%%%%%%%%%%%%%%%%%%%%%%%%%
\subsubsection*{Triplet adjoint in exit space}
%%%%%%%%%%%%%%%%%%%%%%%%%%%%%%%%%%%%%%%%%%%%%%%%%%%%%%%%%%%%%%
%%%%%%%%%%%%%%%%%%%%%%%%%%%%%%%%%%%%%%%%%%%%%%%%%%%%%%%%%%%%%%
%%%%%%%%%%%%%%%%%%%%%%%%%%%%%%%%%%%%%%%%%%%%%%%%%%%%%%%%%%%%%%
%%%%%%%%%%%%%%%%%%%%%%%%%%%%%%%%%%%%%%%%%%%%%%%%%%%%%%%%%%%%%%
In the next lemma
(\cf \cite[Lemmas~5.1, 5.2]{Kurasov09},
\cite[Theorem~3.1]{Dijksma05}) we characterize
the operator
\[
\Amax\dfn\Lmax\cap\cH^2
\]
by considering $\Amax$ as an
extension of the operator $A_0$ in $\cH$ defined by
\[
A_0\dfn L_m\hsum l\,,\quad
l\dfn\{(k,k^\prime)\in\fK^2\vrt
d(k^\prime)=Z_dd(k) \}\,.
\]
The matrix $Z_d$ is the matrix direct sum of
$d$ diagonal $m\times m$ matrices
\[
Z\dfn\diag\{z_1,\ldots,z_m\}\,.
\]
The operator $A_0$ is closed in $\cH$ with
$\rho(A_0)=\rho(L)\cap\rho(l)$,
$\rho(l)=\rho(Z_d)=\bbC\setm\cZ$.
\begin{lem}\label{lem:Amax}
Let $z\in\rho(A_0)$ and
$\wtG_z(\cdot)\dfn p(L)G_z(\cdot)
(=G_z(\cC \cdot))$. Then:
\begin{SL}
\item[1)]
$\Amax=A_0\hsum\whfN_z(\Amax)$ with
$\fN_z(\Amax)=\fN_z(\Lmax)$.
\item[2)]
$\Amax=A_0\hsum\{(\wtG_z(c),z\wtG_z(c)+\kmin(c))\vrt
c\in\bbC^d \}$.
\end{SL}
\end{lem}
\begin{proof}
1) It suffices to verify
$A_0\subseteq\Lmax$. But this follows
from $f+k-g_z(c)\in\fH_{2-m}$ for
$f\in\fH_{m+2}$, $k\in\fK$, $c=c(k)$, where
\begin{equation}
c(k)=(c_s(k))\,,\quad
c_s(k)\dfn\sum_{j}d_{\alpha(s,j)}(k)\,.
\label{eq:ckdef}
\end{equation}

To see $\fN_z(\Amax)=\fN_z(\Lmax)$, first note
that $\fN_z(\Amax)=\fN_z(\Lmax)\cap\cH$.
Now $g_z(\bbC^d)\subseteq\cH$ because
\begin{equation}
\wtP(z)^{-1}(L-z)^{-1}=
\wtP(L)^{-1}(L-z)^{-1}+
\sum_{j}\frac{b_j}{z-z_j}
(L-z_j)^{-1}
\label{eq:xc}
\end{equation}
by \eqref{eq:sum}; \cf \cite[Eq.~(4.10)]{Kurasov09}.

2) Straightforward from \eqref{eq:xc}.
\end{proof}
By Lemma~\ref{lem:Amax} the boundary form of $\Amax$
reads
\begin{align*}
&\braket{f_0+\wtG_z(c),\Amax(f^\prime_0+
\wtG_z(c^\prime))}_\cH-
\braket{\Amax(f_0+\wtG_z(c)),f^\prime_0+
\wtG_z(c^\prime)}_\cH
\\
&
=\braket{d(k),(\cG_Z-\cG^*_Z)d(k^\prime)}_{\bbC^{md}}
+\braket{\wtGm_0(f_0+\wtG_z(c)),
\wtGm_1(f^\prime_0+\wtG_z(c^\prime))}_{\bbC^d}
\\
&-\braket{\wtGm_1(f_0+\wtG_z(c)),
\wtGm_0(f^\prime_0+\wtG_z(c^\prime))}_{\bbC^d}
\end{align*}
for $f_0=f+k$, $f^\prime_0=f^\prime+k^\prime$
from $\dom A_0$ and for
$c,c^\prime\in\bbC^d$.

The matrix
\[
\cG_Z\dfn\cG Z_d
\]
(with $\cG^*_Z$ its adjoint),
the boundary operator
\begin{equation}
\begin{split}
&\wtGm=(\wtGm_0,\wtGm_1)\co\dom\Amax\lto\bbC^{2d}\,,
\\
&\wtGm_0(f+k+\wtG_z(c))\dfn c\,,
\\
&\wtGm_1(f+k+\wtG_z(c))\dfn\cC^*\braket{\vp,f}
+\wtR(z)c-\cG^*_b d(k)
\end{split}
\label{eq:wtGm}
\end{equation}
and the matrix valued Nevanlinna function
(recall \eg \cite{Derkach17})
\[
\wtR(z)\dfn\cC^*R(z)\cC\,,\quad
R(z)=(R_{ss^\prime}(z)\dfn
\braket{\vp_s,G_{s^\prime}(z) })\,.
\]
\begin{rems}\label{rem:rem}
\begin{SL}
\item[1.]
The function $\wtR$ extends to $z\in\rho(L)$.
\item[2.]
Since $\fH_{m+2}\subseteq\fH_m$ densely
and continuously, and $\braket{\vp_s,\cdot}$
is bounded on $\fH_{m+2}$, $\braket{\vp_s,\cdot}$
has a continuation to $\fH_m$, which we denote
by the same symbol in $R(z)$.
In \cite[Definition~3.1.2]{Albeverio00}
$R(\img)$ is related to an admissible matrix
for functionals of class $\fH_{-2}\setm\fH_{-1}$,
see also \cite{Hassi09-b}. The reason behind all
this is that our analysis in
$\fK$ (or $\cH$) can be transferred by scaling
to the subspace of $\fH_0$
generated by $\whg_s(z_j)\dfn(L-z_j)^{-1}\whvp_s$, with
$\whvp_s\dfn P(L)^{-1/2}\vp_s\in\fH_{-2}\setm\fH_{-1}$;
see Appendix~\ref{sec:ibt} for details. In this way
$R_{ss^\prime}(z)=
\braket{\whvp_s,\whg_{s^\prime}(z) }$,
where $\braket{\whvp_s,\cdot}$ is a continuation
to $\fH_0$.
\item[3.]
With the operator
$\Gamma=(\Gamma_0,\Gamma_1)
\co\dom\Lmin^*\lto\bbC^{2d}$ defined by
\begin{align*}
\Gamma_0(f+G_z(c))\dfn&c\,,\quad f\in\fH_{m+2}\,,
\quad c\in\bbC^d\,,
\\
\Gamma_1(f+G_z(c))\dfn&\braket{\vp,f}+R(z)c
\end{align*}
the triple $\Pi_\Gamma=(\bbC^d,\Gamma_0,\Gamma_1)$
is a boundary triple for $\Lmin^*$ with
the $\gamma$-field
$z\mapsto\gamma_\Gamma(z)=G_z(\cdot)$
and the Weyl function $M_\Gamma=R$.
\end{SL}
\end{rems}
At this point one makes an assumption
in \cite{Kurasov09}
that the matrix $\cG_Z$ is Hermitian or,
equivalently:
\begin{prop}\label{prop:Equiv}
The three statements are equivalent:
\begin{SL}
\item[$(i)$]
The matrix $\cG_Z$ is Hermitian.
\item[$(ii)$]
The Gram matrix $\cG$ is diagonal in $j$,
and $\cZ\subseteq\bbR\cap\rho(L)$.
\item[$(iii)$]
The Nevanlinna function $R$ satisfies
$R(z_j)\equiv\cR=\cR^*$ for all
$j$.
\end{SL}
\end{prop}
The reason is: For an Hermitian $\cG_Z$
the adjoint in $\cH$ of $\Amax$,
\ie the operator
\[
\Amin\dfn\Amax^*
\]
is closed densely defined symmetric in
$\cH$, and has defect numbers $(d,d)$. Subsequently,
one applies standard theory of extensions
of symmetric operators, and then characterizes
a generalized resolvent of $\Lmin$
associated with a self-adjoint operator $A_0$
in $\cH$.

In order to prove Proposition~\ref{prop:Equiv}
the easiest way is to use the matrix notation
\[
\cM=
(\cM_{s i^\prime}\dfn
R_{s s^\prime}(z_{j^\prime}))
\]
and to observe that
\begin{equation}
(\ol{z_j}-z_{j^\prime})\cG_{ii^\prime}=
(\cG^*_Z-\cG_Z)_{ii^\prime}=
(\cM^*)_{is^\prime}-
\cM_{si^\prime}\,.
\label{eq:cM}
\end{equation}

For later reference
note that $\cM$ also appears in
\begin{equation}
\begin{split}
&\braket{d(k),(\cG_Z-\cG^*_Z)d(k^\prime)}_{\bbC^{md}}
\\
&=
\braket{c(k),\cM d(k^\prime) }_{\bbC^d}-
\braket{\cM d(k),c(k^\prime) }_{\bbC^d}
\end{split}
\label{eq:cMck}
\end{equation}
with $c(k)$ as in \eqref{eq:ckdef}.
Thus $\cM d(k)=\cR c(k)$ if $\cG^*_Z=\cG_Z$.

In our approach we do not assume that
$\cG_Z$ is necessarily Hermitian. In this case
define the operator
$\Amax^\prime$ in $\cH$ by
\begin{equation}
\Amax^\prime\dfn
A^*_0\hsum\{(\wtG_z(c),z\wtG_z(c)+\kmin(c))\vrt
c\in\bbC^d \}
\label{eq:Amaxp}
\end{equation}
and let $\Amin$ be as previously.
By standard procedure one verifies that
the adjoint $A^*_0$ in $\cH$ of $A_0$
is given by
\[
A^*_0=L_m\hsum l^*\,,\quad
l^*=\{(k,k^\prime)\in\fK^2\vrt
d(k^\prime)=\cG^{-1}\cG^*_Zd(k) \}
\]
where $l^*$ is the adjoint in $\fK$ of $l$.
Moreover
\begin{lem}\label{lem:Amin}
\begin{SL}
\item[1)]
Consider $\wtGm$ as an operator
from a $\whJ_\cH$-space to a $\whJ_{\bbC^d}$-space,
with domain $\Amax$, and let $\wtGm^+$
be its Krein space adjoint. Then
the operator $(\wtGm^+)^{-1}=(\wtGm_0,\wtGm_1)$
but now with domain $\Amax^\prime$.
\item[2)]
In particular,
$\Amin$ is the domain restriction
to $\ker\wtGm$ of $\Amax^\prime$.
\end{SL}
\end{lem}
Instead of requiring $\cG_Z$ to be Hermitian,
in the next section we consider subspaces from
the scale $\fKmin\subseteq\fK$ such that
\eqref{eq:cMck} vanishes.
Without making any additional assumptions,
such subspaces always exist
if $m\geq2$; notice that $k\in\fKmin$ satisfies
$c(k)=0$.
If $m=1$ then $\fKmin=\fK$ and
\eqref{eq:cMck} vanishes iff $z_1\in\bbR\cap\rho(L)$,
\cf Proposition~\ref{prop:Equiv}; hence
$L$ should be semibounded.
%%%%%%%%%%%%%%%%%%%%%%%%%%%%%%%%%%%%%%%%%%%%%%%%%%%%%%%%%%%%%%
%%%%%%%%%%%%%%%%%%%%%%%%%%%%%%%%%%%%%%%%%%%%%%%%%%%%%%%%%%%%%%
%%%%%%%%%%%%%%%%%%%%%%%%%%%%%%%%%%%%%%%%%%%%%%%%%%%%%%%%%%%%%%
%%%%%%%%%%%%%%%%%%%%%%%%%%%%%%%%%%%%%%%%%%%%%%%%%%%%%%%%%%%%%%
\section{Generalized resolvent}\label{sec:gen}
%%%%%%%%%%%%%%%%%%%%%%%%%%%%%%%%%%%%%%%%%%%%%%%%%%%%%%%%%%%%%%
%%%%%%%%%%%%%%%%%%%%%%%%%%%%%%%%%%%%%%%%%%%%%%%%%%%%%%%%%%%%%%
%%%%%%%%%%%%%%%%%%%%%%%%%%%%%%%%%%%%%%%%%%%%%%%%%%%%%%%%%%%%%%
%%%%%%%%%%%%%%%%%%%%%%%%%%%%%%%%%%%%%%%%%%%%%%%%%%%%%%%%%%%%%%
%%%%%%%%%%%%%%%%%%%%%%%%%%%%%%%%%%%%%%%%%%%%%%%%%%%%%%%%%%%%%%
%%%%%%%%%%%%%%%%%%%%%%%%%%%%%%%%%%%%%%%%%%%%%%%%%%%%%%%%%%%%%%
%%%%%%%%%%%%%%%%%%%%%%%%%%%%%%%%%%%%%%%%%%%%%%%%%%%%%%%%%%%%%%
%%%%%%%%%%%%%%%%%%%%%%%%%%%%%%%%%%%%%%%%%%%%%%%%%%%%%%%%%%%%%%
\subsubsection*{Maximal exit subspace}
%%%%%%%%%%%%%%%%%%%%%%%%%%%%%%%%%%%%%%%%%%%%%%%%%%%%%%%%%%%%%%
%%%%%%%%%%%%%%%%%%%%%%%%%%%%%%%%%%%%%%%%%%%%%%%%%%%%%%%%%%%%%%
%%%%%%%%%%%%%%%%%%%%%%%%%%%%%%%%%%%%%%%%%%%%%%%%%%%%%%%%%%%%%%
%%%%%%%%%%%%%%%%%%%%%%%%%%%%%%%%%%%%%%%%%%%%%%%%%%%%%%%%%%%%%%
Fix $m\geq2$ and
define a subspace $\fKmax$ of $\fK$ by
\[
\fKmax\dfn
\{k\in\fK\vrt
\Im\braket{c(k),\cM d(k)}_{\bbC^d}=0 \}\,.
\]
In order to characterize $\fKmax$
it is convenient to interpret $\fKmax$
as a neutral subspace of
a $W$-space \cites{Gohberg05,Azizov89}
as follows. Define an Hermitian matrix
\[
W\dfn\img(\cG_Z-\cG^*_Z)\,.
\]
Then $(\bbC^{md},[\cdot,\cdot])$ is a $W$-space
with an indefinite inner product
\[
[\xi,\xi^\prime]\dfn
\braket{\xi,W\xi^\prime}_{\bbC^{md}}\,,
\quad
\xi,\xi^\prime\in\bbC^{md}\,.
\]
By definition,
a neutral subspace consists of those
$\xi\in\bbC^{md}$ such that $[\xi,\xi]=0$;
it is maximal if it is not contained properly in
a neutral subspace.
Using $\fK\leftrightarrow\bbC^{md}$
and \eqref{eq:cMck},
by polarization therefore $d(\fKmax)$
is maximal neutral.
By \cite[Theorem~2.3.4]{Gohberg05}, the dimension,
$d^\prime$, of $\fKmax$ satisfies
$d^\prime\leq\min\{\pi(W),\nu(W)\}+\dim\ker W$,
where the number of positive (resp. negative) eigenvalues
of $W$, counting multiplicities, is denoted by
$\pi(W)$ (resp. $\nu(W)$).
\begin{prop}\label{prop:K1}
One has the direct sum decomposition
\[
d(\fKmax)=d(\fK\cap\fH_{2-m})\dsum
\ker W
\]
where
$\fK\cap\fH_{2-m}$, of dimension $(m-1)d$,
is the set of those $k\in\fK$ such that $c(k)=0$.

In particular, $\fKmax=\fK\cap\fH_{2-m}$
if $d=1$ and $W\neq0$.
\end{prop}
\begin{proof}
\textit{Step 1.}
For $k\in\fK$
\[
Lk=lk+\sum_{s}
c_s(k)\vp_s\,.
\]
If $k\in\fK\cap\fH_{1-m}$ then
$Lk\in\fH_{-1-m}$. Since $lk\in\fK$
and $\vp_s\in\fH_{-m-2}\setm\fH_{-m-1}$
we have $c(k)=0$, \ie
\[
\fK\cap\fH_{1-m}\subseteq\fK_*\dfn
\{k\in\fK\vrt c(k)=0 \}\,.
\]
Since $\fK_*\subseteq\fKmax\subseteq\fK$, this shows
\[
\fK\cap\fH_{1-m}=\fK_*\cap\fH_{1-m}=
\fKmax\cap\fH_{1-m}\,.
\]
On the other hand
\[
\fK_*=\vee
\{(L-z_m)^{-1}g_i\vrt j\leq m-1\}
\]
\ie $\fK_*\subseteq\fH_{2-m}$. Therefore
\[
\fK_*=\fK\cap\fH_{2-m}=\fKmax\cap\fH_{2-m}\,.
\]

\textit{Step 2.}
To show the direct sum decomposition of $\fKmax$,
consider
\[
\fK_0\dfn\fK_*\dsum(\fKmax\setm\fK_*)\,.
\]
Clearly $\fK_0\subseteq\fKmax$. For the converse
use that $\fKmax\cap\fK^\bot_*\subseteq
\fKmax\setm\fK_*$. Thus
$\fK_0=\fKmax$.

Since
\[
\ker W\subseteq d(\fKmax)=\{\xi\vrt
[\xi,\xi]=0 \}
\]
and
$(\forall(\xi_i)\in\ker W)$
$\sum_j\xi_i\neq0$, we have
$\ker W\subseteq d(\fKmax\setm\fK_*)$.
For the converse, we apply the dimension
argument. We have
\[
\pi(W)+\nu(W)+d_0=md\,,\quad
d_0\dfn\dim\ker W
\]
and
\[
d^\prime\leq d_*+d_0\,,\quad
d_*\dfn\min\{\pi(W),\nu(W)\}\,.
\]
Since
\[
(m-1)d+d_0\leq d^\prime\leq
d_*+d_0
\]
it therefore suffices to show $(m-1)d=d_*$.

Suppose $\pi(W)\leq\nu(W)$. Then
\[
d_*=\pi(W)\leq md-(d_0+\pi(W))\,.
\]
Now $d_0+\pi(W)\geq d^\prime\geq d$ (the last
$\geq$ uses $\fKmin\subseteq\fKmax$), so
$d_*\leq md-d$.
The case $\pi(W)\geq\nu(W)$ is treated
analogously.

\textit{Step 3.}
If $d=1$ then $d^\prime\geq \dim\fK_*=m-1$.
Since moreover $W\neq0$, $d^\prime\leq m-1$.
\end{proof}
\begin{rems}
\begin{SL}
\item[1.]
$c(k)=0$ for $k\in\fKmin$ is seen directly from
$d(\fKmin)=\ran\hb$ and $\sum_jb_j=0$,
recall \eqref{eq:sum}, \eqref{eq:hb}.
Moreover, by \eqref{eq:cM}
\[
\{0\}=\ran\hb\cap\ker\cM=\ran\hb\cap\ker W\,.
\]
If \eg $m=2$ then
$\fK_*=\fKmin$, \ie
$d(\fKmax)=d(\fKmin)\dsum\ker W$.
\item[2.]
In the course of proving Proposition~\ref{prop:K1}
we have on the way established that
\[
\fK\cap\fH_{1-m}=\fK\cap\fH_{2-m}\,.
\]
\end{SL}
\end{rems}
%%%%%%%%%%%%%%%%%%%%%%%%%%%%%%%%%%%%%%%%%%%%%%%%%%%%%%%%%%%%%%
%%%%%%%%%%%%%%%%%%%%%%%%%%%%%%%%%%%%%%%%%%%%%%%%%%%%%%%%%%%%%%
%%%%%%%%%%%%%%%%%%%%%%%%%%%%%%%%%%%%%%%%%%%%%%%%%%%%%%%%%%%%%%
%%%%%%%%%%%%%%%%%%%%%%%%%%%%%%%%%%%%%%%%%%%%%%%%%%%%%%%%%%%%%%
\subsubsection*{Triplet adjoint in exit subspace}
%%%%%%%%%%%%%%%%%%%%%%%%%%%%%%%%%%%%%%%%%%%%%%%%%%%%%%%%%%%%%%
%%%%%%%%%%%%%%%%%%%%%%%%%%%%%%%%%%%%%%%%%%%%%%%%%%%%%%%%%%%%%%
%%%%%%%%%%%%%%%%%%%%%%%%%%%%%%%%%%%%%%%%%%%%%%%%%%%%%%%%%%%%%%
%%%%%%%%%%%%%%%%%%%%%%%%%%%%%%%%%%%%%%%%%%%%%%%%%%%%%%%%%%%%%%
Let $\fK_1=(\fK_1,\braket{\cdot,\cdot}_{-m})$
be an arbitrary subspace from the scale
$\fKmin\subseteq\fKmax$. The
orthogonal complement in $\fK$ of $\fK_1$ is denoted by
\[
\cH_\bot\dfn\fK\om\fK_1\,.
\]
The corresponding Hilbert subspace of $\cH$ is defined by
\[
\cH_1\dfn(\fH_m\dsum\fK_1,\braket{\cdot,\cdot}_\cH)\,.
\]
Notice that the orthogonal complement in
$\cH$ of $\cH_1$ is $\cH_\bot$.

Define a relation $\Bmax$ in $\cH$ by
\[
\Bmax\dfn\Amax\vrt_{\cH_1}\hsum
(\{0\}\times\cH_\bot)\,.
\]
By Lemma~\ref{lem:Amax} and \eqref{eq:Amaxp}
\[
\ran((\Amax-\Amax^\prime)\vrt_{\cH_1})
\subseteq\cH_\bot
\]
so the proof of the next lemma is accomplished
by standard computation.
\begin{lem}\label{lem:Bmax}
Let $\Bmin\dfn\Bmax^*$ be the adjoint in $\cH$
of $\Bmax$. Then
\[
\Bmin=\Amin\vrt_{\cH_1}\hsum
(\{0\}\times\cH_\bot)
\]
is a closed symmetric relation in $\cH$ with
defect numbers $(d,d)$. The adjoint in $\cH$
is given by
\[
\Bmin^*=\Bmax=B_0\hsum\whfN_z(\Bmax)\,,\quad
z\in\rho(B_0)\,.
\]
The self-adjoint relation $B_0$ in $\cH$ admits
a canonical form
\[
B_0\dfn \mB_0\hop
(\{0\}\times\cH_\bot)\,,\quad
\mB_0\dfn L_m\hsum l_1
\]
where
(with $P_{\fK_1}$ an orthogonal projection
in $\fK$ onto $\fK_1$)
\[
l_1\dfn P_{\fK_1}l\vrt_{\fK_1}
\]
is the self-adjoint
operator in $\fK_1$.
The defect subspace
\[
\fN_z(\Bmax)=
[(L-z)^{-1}-(l_1-z)^{-1}](\fKmin)\,.
\]

Moreover, if $\mBmin$ denotes the operator part of
$\Bmin$, \ie
\[
\mBmin=P_{\cH_1}\Amin\vrt_{\cH_1}
\]
($P_{\cH_1}$ is an orthogonal projection in
$\cH$ onto $\cH_1$), then $\mBmin$ is a closed
densely defined  symmetric operator in $\cH_1$,
with equal defect numbers
$(d,d)$. The adjoint in $\cH_1$,
$\mBmax\dfn\mBmin^*$, is the operator part of
$\Bmax$, \ie
\[
\mBmax=\mB_0\hsum\whfN_z(\mBmax)\,,\quad
\fN_z(\mBmax)=\fN_z(\Bmax)\,,\quad
z\in\rho(B_0)\,.
\]
\end{lem}
The boundary value space of $\Bmax$
in $\cH$ is therefore
completely determined by the boundary value space of
$\mBmax$ in $\cH_1$.
\begin{thm}\label{thm:mBmax}
\begin{SL}
\item[1)]
Consider $\wtGm$ as an operator
$\Amax\lto\bbC^{2d}$.
With an operator
\[
\mGm\dfn(\wtGm\vrt_{\cH_1\times\cH})\hsum
((\{0\}\times\cH_\bot)\times\{0\})\dfn
(\mGm_0,\mGm_1)
\]
from $\Bmax$ to $\bbC^{2d}$, the triple
$\Pi_{\mGm}=
(\bbC^d,\mGm_0,\mGm_1)$ is a boundary triple
for $\Bmax$ with the $\gamma$-field $\gamma_{\mGm}$
and the Weyl function $M_{\mGm}$ given on $\bbC^d$ by
\begin{align*}
\gamma_{\mGm}(z)=&
[(L-z)^{-1}-(l_1-z)^{-1}]\kmin(\cdot)\,,
\\
M_{\mGm}(z)=&\wtR(z)+\mQ(z)\,,\quad
z\in\rho(\mB_0)
\end{align*}
where the Nevanlinna function
\[
\mQ(z)\dfn(\braket{\ff_s,(l_1-z)^{-1}
\ff_{s^\prime}}_{-m})\,,\quad
z\in\rho(l_1)\,.
\]
\item[2)]
With an operator
\[
\mGm\dfn\wtGm\vrt_{\cH_1}\dfn(\mGm_0,\mGm_1)
\]
from $\dom\mBmax$ to $\bbC^{2d}$, the triple
$\Pi_{\mGm}=
(\bbC^d,\mGm_0,\mGm_1)$ is a boundary triple
for $\mBmax$ with the $\gamma$-field
$\gamma_{\mGm}$ and the Weyl function $M_{\mGm}$.
\item[3)]
Consider a closed symmetric operator
\[
\ml_1\dfn l_1\vrt_{\fK_1\om\fKmin}
\]
in $\fK_1$, with defect numbers $(d,d)$;
the adjoint in $\fK_1$ is characterized by
\[
\begin{split}
&\ml^*_1=l_1\hsum(\{0\}\times\fKmin)\,,
\\
&\fN_z(\ml^*_1)=(l_1-z)^{-1}(\fKmin)\,,\quad
z\in\rho(l_1)\,.
\end{split}
\]
Then the triple $\Pi_{\mmGm}=
(\bbC^d,\mmGm_0,\mmGm_1)$, where
\[
\begin{split}
&\mmGm=(\mmGm_0,\mmGm_1)\co
\ml^*_1\lto\bbC^{2d}\,,
\\
&\mmGm_0(k,l_1k+\kmin(\chi))\dfn \chi\,,
\\
&\mmGm_1(k,l_1k+\kmin(\chi))\dfn
-\cG^*_b d(k)
\end{split}
\]
is a boundary triple for $\ml^*_1$ with
the $\gamma$-field $\gamma_{\mmGm}$
and the Weyl function $M_{\mmGm}$
given on $\bbC^d$ by
\[
\gamma_{\mmGm}(z)=-(l_1-z)^{-1}\kmin(\cdot) \,,\quad
M_{\mmGm}(z)=\mQ(z)\,,\quad
z\in\rho(l_1)\,.
\]
\end{SL}
\end{thm}
\begin{proof}
1) $\mGm$ can be given the form
\[
\mGm=(\wtGm\cap\fM)\hsum\fN
\]
where the closed relations
$\fM$ and $\fN$ from
$\cH^2$ to $\bbC^{2d}$ are defined by
\[
\fM\dfn(\cH_1\times\cH)\times\bbC^{2d}\,,\quad
\fN\dfn(\{0\}\times\cH_\bot)\times\{0\}\,.
\]
In order that $\Pi_{\mGm}$ should be
the boundary triple for $\Bmax$
it is necessary and sufficient that the operator
$\mGm$ should be unitary from a
$\whJ_\cH$-space to a $\whJ_{\bbC^d}$-space.
Since the Krein space adjoints
$\fN^+=\fM^{-1}$ and $\fM^+=\fN^{-1}$,
the Krein space adjoint of $\mGm$ is given by
\[
\mGm^+=(\wtGm\cap\fM)^+\cap\fN^+=
((\ol{(\wtGm^+)^{-1}\hsum\fN})\cap\fM)^{-1}\,.
\]
By Lemma~\ref{lem:Amin}
\[
(\wtGm^+)^{-1}\hsum\fN=\wtGm\hsum\fN\,.
\]
Because the relation $\wtGm\hsum\fN$ is closed
and $\fN\subseteq\fM$, one therefore gets
$\mGm^+=\mGm^{-1}$ as required.

The computation of $\gamma_{\mGm}$,
$M_{\mGm}$ is standard by applying Lemma~\ref{lem:Bmax}.
Note moreover that $\mQ(z)$ as stated follows
from
\[
\mQ(z)=\cG^*_bd((l_1-z)^{-1}\kmin(\cdot))\,.
\]

2)
With $\wtGm\co\Amax\lto\bbC^{2d}$ as previously,
one needs to verify that the operator
\[
\mGm\dfn((\wtGm\cap\fM)\hsum\fN)\vrt_{\cH^2_1}
\]
is unitary from a
$\whJ_{\cH_1}$-space to a $\whJ_{\bbC^d}$-space;
observe that the above $\mGm$ coincides with
$\wtGm\vrt_{\cH_1}$ if considered as
an operator $\dom\mBmax\lto\bbC^{2d}$.
Thus, the inverse of the Krein space adjoint
\begin{align*}
(\mGm^+)^{-1}=&
(\ol{\mGm\hsum(\cH^2_\bot\times\{0\} ) })
\vrt_{\cH^2_1}
\\
=&(\mGm\hsum(\cH^2_\bot\times\{0\} ) )
\vrt_{\cH^2_1}
\\
=&((\wtGm\cap\fM)\hsum\fN)
\vrt_{\cH^2_1}=\mGm
\end{align*}
as required.

3) Routine computation.
\end{proof}
\begin{rem}
With $\mGm$ as in 1),
$\mmGm\supsetneq\mGm\vrt_{\fK^2_1}$; in fact
\[
\mmGm=\mGm\vrt_{\fK^2_1}\hsum
\{((0,\kmin(\chi)),(\chi,0))\vrt
\chi\in\bbC^d \}\,.
\]
\end{rem}
The subsequent
results are presented only for $\mBmax$.
%%%%%%%%%%%%%%%%%%%%%%%%%%%%%%%%%%%%%%%%%%%%%%%%%%%%%%%%%%%%%%
%%%%%%%%%%%%%%%%%%%%%%%%%%%%%%%%%%%%%%%%%%%%%%%%%%%%%%%%%%%%%%
%%%%%%%%%%%%%%%%%%%%%%%%%%%%%%%%%%%%%%%%%%%%%%%%%%%%%%%%%%%%%%
%%%%%%%%%%%%%%%%%%%%%%%%%%%%%%%%%%%%%%%%%%%%%%%%%%%%%%%%%%%%%%
\subsubsection*{Generalized resolvent}
%%%%%%%%%%%%%%%%%%%%%%%%%%%%%%%%%%%%%%%%%%%%%%%%%%%%%%%%%%%%%%
%%%%%%%%%%%%%%%%%%%%%%%%%%%%%%%%%%%%%%%%%%%%%%%%%%%%%%%%%%%%%%
%%%%%%%%%%%%%%%%%%%%%%%%%%%%%%%%%%%%%%%%%%%%%%%%%%%%%%%%%%%%%%
%%%%%%%%%%%%%%%%%%%%%%%%%%%%%%%%%%%%%%%%%%%%%%%%%%%%%%%%%%%%%%
A closed operator $\mB\in\mrm{Ext}(\mBmin)$
is in bijective correspondence with
a closed relation $\Theta$ in $\bbC^d$
via $\mB=\mB_\Theta\dfn\mGm^{-1}(\Theta)$, \ie
\[
\mB_\Theta\subseteq \mBmax\,,\quad
\dom \mB_\Theta\dfn\{f\in\dom\mBmax\vrt
(\mGm_0f,\mGm_1f)\in\Theta \}\,.
\]
Thus $\mBmin=\mB_{\Theta=\{0\}}$
and $\mB_0=\mB_{\Theta=\{0\}\times\bbC^d}$.

In the next theorem
$U\co f+k\mapsto(f,k)$
is a unitary operator from a Hilbert space
$\cH_1$ to an (external) Hilbert sum
$\fH_m\op\fK_1$.
\begin{thm}\label{thm:gen}
\begin{SL}
\item[1)]
Let $\Theta$ be a closed relation in $\bbC^d$.
The resolvent of a closed operator
$\mB_\Theta$
is given by
\[
U(\mB_\Theta-z)^{-1}U^*=
\begin{pmatrix}
\mR^{11}_\Theta(z) & \mR^{12}_\Theta(z) \\
\mR^{21}_\Theta(z) & \mR^{22}_\Theta(z)
\end{pmatrix}\co
\begin{matrix}
\fH_m \\ \op \\ \fK_1
\end{matrix}
\lto
\begin{matrix}
\fH_m \\ \op \\ \fK_1
\end{matrix}
\]
for $z\in\rho(\mB_\Theta)\cap\rho(\mB_0)$;
the entries
{\small
\begin{align*}
\mR^{11}_\Theta(z)\dfn&
(L-z)^{-1}+G_z
\bigl(\cC(\Theta-M_{\mGm}(z))^{-1}\cC^*
\braket{\vp,(L-z)^{-1}\cdot}\bigr)
\,,
\\
\mR^{12}_\Theta(z)\dfn&
-G_z
\bigl(\cC(\Theta-M_{\mGm}(z))^{-1}\cG^*_b
d((l_1-z)^{-1}\cdot) \bigr)
\,,
\\
\mR^{21}_\Theta(z)\dfn&
-(l_1-z)^{-1}\kmin
\bigl((\Theta-M_{\mGm}(z))^{-1}\cC^*
\braket{\vp,(L-z)^{-1}\cdot} \bigr)\,,
\\
\mR^{22}_\Theta(z)\dfn&
(l_1-z)^{-1}+
(l_1-z)^{-1}\kmin
\bigl((\Theta-M_{\mGm}(z))^{-1}
\cG^*_bd((l_1-z)^{-1}\cdot) \bigr)\,.
\end{align*}}
\item[2)]
Let $\Theta$ be a self-adjoint relation
in $\bbC^d$.
To a generalized resolvent
$P_{\fH_m}(\mB_\Theta-z)^{-1}\vrt_{\fH_m}$
there corresponds, via the generalized Krein--Naimark
resolvent formula, a Nevanlinna family
\[
\tau\co
z\mapsto(\cC^{-1})^*(\mQ(z)-\Theta)\cC^{-1}\,.
\]
\end{SL}
\end{thm}
\begin{proof}
1) This part is due to
Lemma~\ref{lem:Bmax}, Theorem~\ref{thm:mBmax},
and the Krein--Naimark resolvent formula.

2) We need to verify that
the only solution $\tau(z)$ to
\[
\mR^{11}_\Theta(z)=(L-z)^{-1}-
\gamma_\Gamma(z)(\tau(z)+M_\Gamma(z))
\gamma_\Gamma(\ol{z})^*
\]
is as stated; see Remark~\ref{rem:rem}-3.
The above equation reads
\begin{align*}
&G_z
(\cC(\Theta-M_{\mGm}(z))^{-1}\cC^*
\braket{\vp,(L-z)^{-1}\cdot})
\\
&=-G_z
((\tau(z)+M_{\Gamma}(z))^{-1}
\braket{\vp,(L-z)^{-1}\cdot})
\end{align*}
and then
\[
\cC(\Theta-M_{\mGm}(z))^{-1}\cC^*
=-(\tau(z)+R(z))^{-1}
\]
and
\[
\cC^{*\,-1}(\Theta-M_{\mGm}(z))\cC^{-1}
=-(\tau(z)+R(z))
\]
with both sides considered as relations.
Subsequently, a relation $\tau(z)$ is
the operatorwise sum of relations
\[
\cC^{*\,-1}M_{\mGm}(z)\cC^{-1}-R(z)=
\cC^{*\,-1}\mQ(z)\cC^{-1}
\]
and $-\cC^{*\,-1}\Theta\cC^{-1}$.
\end{proof}
\begin{rem}\label{rem:rem2}
If $d=1$, $L\geq0$, $\cZ=\cZ\cap(-\infty,0)$,
$\cG$ is diagonal,
then an operator $\mB_\Theta=B_\Theta$ in $\cH$,
and then the resolvent in Theorem~\ref{thm:gen}
is unitarily equivalent to that in
\cite[Theorem~6.1]{Kurasov09}.
The scalar $Q$-function of the symmetric
operator $\ml_1=l\vrt_{\fK\om\fKmin}$ in $\fK$
with defect numbers $(1,1)$ is
\[
\mQ(z)=\braket{b,\cG(Z-z)^{-1}b}_{\bbC^m}\,,
\quad
b\dfn(b_j)\in\bbC^m\,.
\]
To compare with,
the $Q$-function of the symmetric
operator $\ml_1=\{0\}$ in $\fKmin$
(\ie $\fK_1=\fKmin$) is given by
\[
\mQ(z)=\braket{b,\cG b}^2_{\bbC^m}/
\braket{b,\cG(Z-z)b}_{\bbC^m}\,.
\]
Notice that
$(\cG Z-Z^*\cG)b\bot b$, and
$\cG Z\neq Z^*\cG$ is allowed in this case.
\end{rem}
%%%%%%%%%%%%%%%%%%%%%%%%%%%%%%%%%%%%%%%%%%%%%%%%%%%%%%%%%%%%%%
%%%%%%%%%%%%%%%%%%%%%%%%%%%%%%%%%%%%%%%%%%%%%%%%%%%%%%%%%%%%%%
%%%%%%%%%%%%%%%%%%%%%%%%%%%%%%%%%%%%%%%%%%%%%%%%%%%%%%%%%%%%%%
%%%%%%%%%%%%%%%%%%%%%%%%%%%%%%%%%%%%%%%%%%%%%%%%%%%%%%%%%%%%%%
\section{Final remarks}\label{sec:fin}
%%%%%%%%%%%%%%%%%%%%%%%%%%%%%%%%%%%%%%%%%%%%%%%%%%%%%%%%%%%%%%
%%%%%%%%%%%%%%%%%%%%%%%%%%%%%%%%%%%%%%%%%%%%%%%%%%%%%%%%%%%%%%
%%%%%%%%%%%%%%%%%%%%%%%%%%%%%%%%%%%%%%%%%%%%%%%%%%%%%%%%%%%%%%
%%%%%%%%%%%%%%%%%%%%%%%%%%%%%%%%%%%%%%%%%%%%%%%%%%%%%%%%%%%%%%
%%%%%%%%%%%%%%%%%%%%%%%%%%%%%%%%%%%%%%%%%%%%%%%%%%%%%%%%%%%%%%
%%%%%%%%%%%%%%%%%%%%%%%%%%%%%%%%%%%%%%%%%%%%%%%%%%%%%%%%%%%%%%
%%%%%%%%%%%%%%%%%%%%%%%%%%%%%%%%%%%%%%%%%%%%%%%%%%%%%%%%%%%%%%
%%%%%%%%%%%%%%%%%%%%%%%%%%%%%%%%%%%%%%%%%%%%%%%%%%%%%%%%%%%%%%
\subsubsection*{Minimal realization}
%%%%%%%%%%%%%%%%%%%%%%%%%%%%%%%%%%%%%%%%%%%%%%%%%%%%%%%%%%%%%%
%%%%%%%%%%%%%%%%%%%%%%%%%%%%%%%%%%%%%%%%%%%%%%%%%%%%%%%%%%%%%%
%%%%%%%%%%%%%%%%%%%%%%%%%%%%%%%%%%%%%%%%%%%%%%%%%%%%%%%%%%%%%%
%%%%%%%%%%%%%%%%%%%%%%%%%%%%%%%%%%%%%%%%%%%%%%%%%%%%%%%%%%%%%%
The Nevanlinna function $\mQ$ can be
given a standard form (\cf \eqref{eq:KQ})
\[
\mQ(z)=\mQ(\ol{z_0})+(z-\ol{z_0})
\gamma_{\mmGm}(z_0)^*
(I+(z-z_0)(l_1-z)^{-1} )\gamma_{\mmGm}(z_0)
\]
for $z$, $z_0\in\rho(l_1)$. One says in
\cites{Hassi16,Dijksma04a} the pair
$(l_1,\gamma_{\mmGm})$ realizes $\mQ$, and
a realization is minimal if
($\vee$ $\equiv$ closed linear span)
\[
\fK_1=\vee\{
(I+(z-z_0)(l_1-z)^{-1} )\gamma_{\mmGm}(z_0)\chi
\vrt\chi\in\bbC^d\,;\,z\in\rho(l_1) \}\,.
\]
Since the right-hand side
$=\vee\{\fN_z(\ml^*_1)\vrt z\in\rho(l_1) \}$,
this means $\ml_1$ should be simple, \ie
$\ml_1$ should not admit orthogonal decompositions
with a self-adjoint summand.
\begin{prop}\label{prop:min}
$\ml_1$ is simple iff $\fK_1=\fKmin$.
\end{prop}
\begin{proof}
Sufficiency is clear, so we prove necessity.
Consider
$\fK_1$ as a set generated by elements
$\{\fh_\mu\}_{\mu=1}^{d_1}$,
with the corresponding Gram matrix
\[
\cG_1\dfn(\braket{\fh_\mu,\fh_{\mu^\prime} }_{-m})
\]
and $d_1\dfn\dim\fK_1$.
An element $k$ from $\fK_1$ is of
the form
\[
k=k(\eta)\dfn\sum_{\mu}\eta_\mu \fh_\mu\,,
\quad
\eta=(\eta_\mu)=\cG^{-1}_1\cW^*d(k)
\]
with $d(\fK_1)=\cG^{-1}(\ran\cW)$. The
transfer matrix
\[
\cW\dfn(\braket{g_i,\fh_\mu }_{-m})\,.
\]
Subsequently, $l_1$ maps
$k(\eta)$ to $k(\mZ_{d_1}\eta)$, where the
matrix
\begin{equation}
\mZ_{d_1}\dfn\cG^{-1}_1\cW^*Z_d\cG^{-1}\cW
\label{eq:mZd1}
\end{equation}
is similar to an Hermitian matrix:
By recalling that $d(\fK_1)$ is a neutral subspace of a
$W$-space it holds
\begin{equation}
\mZ^*_{d_1}\cG_1=\cG_1\mZ_{d_1}\,.
\label{eq:mZd1-b}
\end{equation}

Let $\{\lambda_\mu\}$ be the eigenvalues
of $\mZ_{d_1}$. By using the spectral decomposition
\[
\mZ_{d_1}=V_1^{-1}\diag\{\lambda_\mu\}V_1\,,
\quad
V_1\dfn U_1\cG^{1/2}_1
\]
with a $d_1\times d_1$ unitary matrix $U_1$,
the eigenspace $\fN_\lambda(\ml_1)$,
for some $\lambda\in\bbR$, consists of
$k(\eta)\in\fK_1\om\fKmin$ such that
$(\forall\mu)$ $(\lambda_\mu-\lambda)(V_1\eta)_\mu=0$.
Now, a simple $\ml_1$ has no eigenvalues,
so that necessarily $\fK_1=\fKmin$.
\end{proof}
\begin{rem}\label{rem:rem3}
With the notation as in Proposition~\ref{prop:min}
$\mQ$ can be given a yet another form
\[
\mQ(z)=
(\hb^*\cW)(\mZ_{d_1}-z)^{-1}
\cG^{-1}_1(\hb^*\cW)^*\,.
\]
\end{rem}
%%%%%%%%%%%%%%%%%%%%%%%%%%%%%%%%%%%%%%%%%%%%%%%%%%%%%%%%%%%%%%
%%%%%%%%%%%%%%%%%%%%%%%%%%%%%%%%%%%%%%%%%%%%%%%%%%%%%%%%%%%%%%
%%%%%%%%%%%%%%%%%%%%%%%%%%%%%%%%%%%%%%%%%%%%%%%%%%%%%%%%%%%%%%
%%%%%%%%%%%%%%%%%%%%%%%%%%%%%%%%%%%%%%%%%%%%%%%%%%%%%%%%%%%%%%
\subsubsection*{Invariance}
%%%%%%%%%%%%%%%%%%%%%%%%%%%%%%%%%%%%%%%%%%%%%%%%%%%%%%%%%%%%%%
%%%%%%%%%%%%%%%%%%%%%%%%%%%%%%%%%%%%%%%%%%%%%%%%%%%%%%%%%%%%%%
%%%%%%%%%%%%%%%%%%%%%%%%%%%%%%%%%%%%%%%%%%%%%%%%%%%%%%%%%%%%%%
%%%%%%%%%%%%%%%%%%%%%%%%%%%%%%%%%%%%%%%%%%%%%%%%%%%%%%%%%%%%%%
In order to emphasize that $\mQ$ is associated to
$(\ml_1,l_1)$ in $\fK_1$, we use the notation
\[
\mQ\equiv\mQ_{\fK_1}\,.
\]
Because the diagonal entry
\cite[Lemma~2.3]{Gallone20}
\[
(\mQ_{\fK_1}(\lambda))_{ss}=
\sup_{k\in\fK_1}
\frac{\abs{\braket{k,\ff_s}_{-m} }^2 }{
\braket{k,(l-\lambda)k}_{-m}}\,,
\quad \lambda<\min\sigma(l_1)
\]
(as usual the spectrum $\sigma(\cdot)\dfn
\bbC\setm\rho(\cdot)$)
and (\eg \cite{Behrndt13})
a scalar valued Nevanlinna function is monotonically
nondecreasing on any interval of $\bbR$ where
it is analytic, one concludes that
(\cf Remark~\ref{rem:rem2})
\begin{equation}
(\mQ_{\fK^\prime_1}(\lambda))_{ss}
\leq (\mQ_{\fK_1}(\lambda))_{ss}\quad
\text{if}\quad
\fK^\prime_1\subseteq\fK_1
\label{eq:x}
\end{equation}
where $\fK^\prime_1$ is some other subspace
from the scale $\fKmin\subseteq\fKmax$.
By the min-max principle
$\min\sigma(l^\prime_1)
\geq\min\sigma(l_1)$,
$l^\prime_1\dfn P_{\fK^\prime_1}l\vrt_{\fK^\prime_1}$.
The preceding considerations imply the following
spectral feature in case $d=1$:
The points
$\{\lambda\vrt \mQ_{\fK_1}(\lambda)=\Theta-\wtR(\lambda)\}$
in $\rho(L)\cap(-\infty,\min\sigma(l_1))$
from the spectrum
of a self-adjoint operator $\mB_\Theta$
in $\cH_1$,
$\abs{\Theta}<\infty$,
are nondecreasing whenever $\fK_1$ gets smaller.

These points need not
increase, as $\fK_1$ gets smaller, if
$\fK_1$ is invariant for $l$, \ie if
$l(\fK_1)\subseteq\fK_1$.
In this case $l=l_1\op (l\vrt_{\cH_\bot})$,
so that a generalized
resolvent $P_{\fK_1}(l-z)^{-1}\vrt_{\fK_1}$
coincides with $(l_1-z)^{-1}$
(and is called orthogonal in
\cite[Chapter~2]{Albeverio00}); hence
$\mQ_{\fK_1}=\mQ_\fK$.
However
\begin{prop}
$\fK_1$ is an invariant subspace for $l$
iff $\fK_1=\fK$.
\end{prop}
\begin{proof}
If $l(\fK_1)\subseteq\fK_1$
or, equivalently, $Z_dd(\fK_1)\subseteq d(\fK_1)$,
then by Halmos \cite[Theorem~3]{Halmos71}
there is an $md\times md$ matrix
$C=(C_{ii^\prime})$
that commutes with $Z_d$ and satisfies
$d(\fK_1)=\ker C$. The commutation criterion
shows that $C$ is diagonal in $j$
(as previously, $i=\alpha(s,j)$,
$i^\prime=\alpha(s^\prime,j^\prime)$).
On the other hand
\[
\ker C\supseteq
\ran\hb=
\vee\{\sum_{j=1}^{m-1}
b_j(e_{sm-m+j}-e_{sm}) \}
\]
where $\{e_{i}\}_{i=1}^{md}$
is a standard basis of $\bbC^{md}$,
so that then
$\sum_jb_jCe_{i}=0$.
But $Ce_i=(C_{1,i},\ldots,C_{md,i})$
and $C_{1,i}=\delta_{j,1}C_{1,i},\ldots,
C_{md,i}=\delta_{j,m}C_{md,1} $, thus
($\forall i$) $b_jCe_i=0$, \ie $C=0$.
\end{proof}
%%%%%%%%%%%%%%%%%%%%%%%%%%%%%%%%%%%%%%%%%%%%%%%%%%%%%%%%%%%%%%
%%%%%%%%%%%%%%%%%%%%%%%%%%%%%%%%%%%%%%%%%%%%%%%%%%%%%%%%%%%%%%
%%%%%%%%%%%%%%%%%%%%%%%%%%%%%%%%%%%%%%%%%%%%%%%%%%%%%%%%%%%%%%
%%%%%%%%%%%%%%%%%%%%%%%%%%%%%%%%%%%%%%%%%%%%%%%%%%%%%%%%%%%%%%
\subsubsection*{Generalized resolvent of $l_1$}
%%%%%%%%%%%%%%%%%%%%%%%%%%%%%%%%%%%%%%%%%%%%%%%%%%%%%%%%%%%%%%
%%%%%%%%%%%%%%%%%%%%%%%%%%%%%%%%%%%%%%%%%%%%%%%%%%%%%%%%%%%%%%
%%%%%%%%%%%%%%%%%%%%%%%%%%%%%%%%%%%%%%%%%%%%%%%%%%%%%%%%%%%%%%
%%%%%%%%%%%%%%%%%%%%%%%%%%%%%%%%%%%%%%%%%%%%%%%%%%%%%%%%%%%%%%
We present a (unique) Nevanlinna
family $\tau_1$ corresponding to
$P_{\fKmin}(l_1-z)^{-1}\vrt_{\fKmin}$ in
terms of $\mQ_{\fK_1}$.
Below we use
\[
\mZ_d=\cGmin^{-1}\cG^*_bZ_d\hb
\]
as in \eqref{eq:mZd1} for $\fK_1=\fKmin$.
\begin{prop}
To a generalized resolvent
$P_{\fKmin}(l_1-z)^{-1}\vrt_{\fKmin}$
there corresponds, via the generalized Krein--Naimark
resolvent formula, a matrix valued Nevanlinna function
\[
\tau_1\co z\mapsto
-\cGmin(\mQ_{\fK_1}(z)^{-1}\cGmin+z)\,.
\]

If $\fK_1=\fKmin$ then
$\tau_1(z)=-\cGmin\mZ_d$ is an Hermitian matrix.
\end{prop}
\begin{proof}
Consider a simple symmetric
operator $\{0\}$ in $\fKmin$. The boundary triple
$(\bbC^d,\Gamma^{\min}_0,\Gamma^{\min}_1)$
for $\fKmin^2$ is given by
($\chi,\chi^\prime\in\bbC^d$)
\begin{align*}
\Gamma^{\min}_0(\kmin(\chi),\kmin(\chi^\prime))
\dfn&\chi\,,
\\
\Gamma^{\min}_1(\kmin(\chi),\kmin(\chi^\prime))
\dfn&\cGmin\chi^\prime\,.
\end{align*}
The $\gamma$-field and the Weyl function
corresponding to the triple are given by
$\kmin(\cdot)$ and $z\cGmin$, respectively.
Subsequently, by using
$\ker\Gamma^{\min}_0=\{0\}\times\fKmin$
and by applying the generalized Krein--Naimark
formula
\[
\mQ_{\fK_1}(z)=-\cGmin
(\tau_1(z)+z\cGmin)^{-1}\cGmin
\]
for $z\in\rho(l_1)$. Now $\mQ_{\fK_1}(z)$,
as the Weyl family corresponding to
$\Pi_{\mmGm}$, is nonsingular (or use that
$l_1$ and $\{0\}\times\fKmin$ are disjoint),
so $\tau_1(z)$ is the matrix as stated.

If $\fK_1=\fKmin$ then use that
\[
\mQ_{\fKmin}(z)=\cGmin(\mZ_d-z)^{-1}
\]
for $z\in\rho(\mZ_d)$, \cf Remark~\ref{rem:rem3},
and recall from \eqref{eq:mZd1-b}
that the matrix $\cGmin\mZ_d$ is Hermitian.
\end{proof}
\begin{rem}
If $d=1$ then the scalar
\[
\mQ_{\fK_1}(z)=
\frac{z\norm{\ff}^{2}_{-m}-\braket{\ff,L\ff}_{-m}}{
z\norm{\ff}^{2}_{-m}+\tau_1(z)}\;
\mQ_{\fKmin}(z)\,.
\]
As already pointed out in Remark~\ref{rem:rem2}
$\braket{b,\cG Z b}_{\bbC^m}$ is real valued;
actually it equals
$\braket{\ff,L\ff}_{-m}$ (with $\ff\equiv\ff_s$ for $d=1$).
Similarly $\braket{b,\cG b}_{\bbC^m}=
\norm{\ff}^2_{-m}$.
The scalar $\tau_1(z)=-\braket{\ff,L\ff}_{-m}$
iff $\fK_1=\fKmin$; otherwise, by using \eqref{eq:x}
$\tau_1(\lambda)\leq-\braket{\ff,L\ff}_{-m}$
for $\lambda<\min\sigma(l_1)$.
\end{rem}
\begin{appendix}
%%%%%%%%%%%%%%%%%%%%%%%%%%%%%%%%%%%%%%%%%%%%%%%%%%%%%%%%%%%%%%
%%%%%%%%%%%%%%%%%%%%%%%%%%%%%%%%%%%%%%%%%%%%%%%%%%%%%%%%%%%%%%
%%%%%%%%%%%%%%%%%%%%%%%%%%%%%%%%%%%%%%%%%%%%%%%%%%%%%%%%%%%%%%
%%%%%%%%%%%%%%%%%%%%%%%%%%%%%%%%%%%%%%%%%%%%%%%%%%%%%%%%%%%%%%
\section{The peak model in the
reference space}
\label{sec:ibt}
%%%%%%%%%%%%%%%%%%%%%%%%%%%%%%%%%%%%%%%%%%%%%%%%%%%%%%%%%%%%%%
%%%%%%%%%%%%%%%%%%%%%%%%%%%%%%%%%%%%%%%%%%%%%%%%%%%%%%%%%%%%%%
%%%%%%%%%%%%%%%%%%%%%%%%%%%%%%%%%%%%%%%%%%%%%%%%%%%%%%%%%%%%%%
%%%%%%%%%%%%%%%%%%%%%%%%%%%%%%%%%%%%%%%%%%%%%%%%%%%%%%%%%%%%%%
We are given a linearly independent system
$\{\whvp_s\in\fH_{-2}\setm\fH_{-1}\}$ and
a closed densely defined
symmetric restriction $\whL_0$ of $L$
subject to the boundary condition $(\forall s)$
$\braket{\whvp_s,u}=0$, $u\in\fH_2$.
The adjoint $\whL^*_0$
in $\fH_0$ of $\whL_0$ extends $L$ to
$\fH_2\dsum\fN_z(\whL^*_0)$, $z\in\rho(L)$,
where the defect subspace
\[
\fN_z(\whL^*_0)=\whg_z(\bbC^d)\,,\quad
\whg_z(c)\dfn\sum_{s}c_s
\whg_s(z)\,,\quad
\whg_s(z)\dfn(L-z)^{-1}\whvp_s
\]
and $c=(c_s)\in\bbC^d$.

Here we analyze the peak model
transformed to $\fH_0$ by using a unitary
operator $P(L)^{-1/2}$ from $\fH_{-n}$ onto
$\fH_{m+n}$; $n\geq0$.
Thus, we view $\whvp_s$ as
\[
\whvp_s\dfn P(L)^{-1/2}\vp_s\,,\quad
\vp_s\in\fH_{-m-2}\setm\fH_{-m-1}
\]
for some fixed integer $m\geq1$.
Then
\[
\whL_0=P(L)^{1/2}\Lmin P(L)^{-1/2}\,,\quad
\whL^*_0=P(L)^{1/2}\Lmin^* P(L)^{-1/2}\,.
\]
The triple
$\Pi_{\whGm}=(\bbC^d,\whGm_0,\whGm_1)$
with (recall Remark~\ref{rem:rem})
\[
\whGm=(\whGm_0,\whGm_1)\dfn
\Gamma P(L)^{-1/2}\,,\quad
\dom\whGm=\dom\whL^*_0
\]
is a boundary triple for $\whL^*_0$ with the
$\gamma$-field $\rho(L)\ni z\mapsto
\gamma_{\whGm}(z)=\whg_z(\cdot)$
and the Weyl function $M_{\whGm}=R$.
In this way $R$ is the
$Q$-function for both $(\Lmin,L_m)$ and
$(\whL_0,L)$.

What we really want to show is that
the scaled boundary operator $\wtGm$ in
\eqref{eq:wtGm} defines an essentially unitary
boundary triple
$\Pi_{\whGm^\prime}=(\bbC^d,\whGm^\prime_0,
\whGm^\prime_1)$ for $\whL^*_0$. This means
the operator
$\whGm^\prime\dfn(\whGm^\prime_0,\whGm^\prime_1)$,
with domain
\[
\whAmax\dfn P(L)^{-1/2}\Amax P(L)^{1/2}
\]
dense in $\whL^*_0$,
is an isometry from a $\whJ_{\fH_0}$-space
to a $\whJ_{\bbC^d}$-space, $\whGm^{\prime\,-1}
\subseteq\whGm^{\prime\,+}$,
and moreover
the closure $\clo\whGm^\prime=\whGm$.
As usual, we consider $\whGm^\prime$ also
as an operator $\dom\whAmax\lto\bbC^{2d}$.
An essentially unitary boundary triple is
a special case of an isometric boundary triple
studied in \cite{Derkach17}.

To make our statement precise,
define a subset $\whfK$ of $\fH_0$ by
\[
\whfK\dfn\vee\{\whg_i\dfn\whg_s(z_j)\}\,,\quad
i=\alpha(s,j)\,.
\]
Note $\whfK\cap\fH_{2m-1}=\{0\}$.
Every $\whk\in\whfK$ is in bijective
correspondence with
$d(\whk)=(d_i(\whk))$ via
the Gram matrix $\cG=(\braket{\whg_i,
\whg_{i^\prime}}_0)$. Let (\cf \eqref{eq:ckdef})
\[
c(\whk)=(c_s(\whk))\in\bbC^d\,,\quad
c_s(\whk)\dfn\sum_{j}d_{i}(\whk)\,.
\]
Then $c(\whfK)=\bbC^d$ and moreover
$c(\whk)=0$ iff $\whk\in\whfK\cap\fH_2$.
Let
\[
\fN\dfn\fH_{2m+2}\dsum
\wtP(L)^{-1}(\fN_z(\whL^*_0)) \,.
\]
Then $\fN$ is dense in $\fH_0$:
$\fH_{2m+2}\subseteq\fN\subseteq\fH_{2m}
\subseteq\fH_2 $.
\begin{thm}
Define an operator
$\whGm^\prime\dfn(\whGm^\prime_0,\whGm^\prime_1)$
where
\begin{align*}
\whGm^\prime_0\co& \fN\dsum\whfK\ni
u+\whk\mapsto c(\whk)\,,
\\
\whGm^\prime_1\co& \fN\dsum\whfK\ni
u+\whk\mapsto \braket{\whvp,u}+\cM d(\whk)\,.
\end{align*}
The triple $\Pi_{\whGm^\prime}=(\bbC^d,\whGm^\prime_0,
\whGm^\prime_1)$ is an isometric boundary triple
for $\whL^*_0$ such that
\[
\clo\whGm^\prime=\whGm\,,\quad
\ran\whGm^\prime_0=\bbC^d\,,\quad
T_0\subsetneq\ol{T_0}=L\,.
\]
Here an essentially self-adjoint operator
\[
T_0\dfn\whL^*_0\vrt_{\ker\whGm^\prime_0}=
L\vrt_{\fN\dsum(\whfK\cap\fH_2)}\,.
\]

The $\gamma$-field and the Weyl function
corresponding to $\Pi_{\whGm^\prime}$ are
$z\mapsto\whg_z(\cdot)$ and $R$, respectively.
\end{thm}
\begin{proof}
It is rather straightforward that
\[
\whAmax(u+\whk)=Lu+
\sum_{i}(Z_dd(\whk))_i\whg_i
\]
and then the boundary form
\begin{align*}
&\braket{u+\whk,\whAmax(u^\prime+\whk^\prime) }_0-
\braket{\whAmax(u+\whk),u^\prime+\whk^\prime }_0
\\
&=
\braket{\whGm^\prime_0(u+\whk),
\whGm^\prime_1(u^\prime+\whk^\prime) }_{\bbC^d}-
\braket{\whGm^\prime_1(u+\whk),
\whGm^\prime_0(u^\prime+\whk^\prime) }_{\bbC^d}
\end{align*}
for $u,u^\prime\in\fN$ and
$\whk,\whk^\prime\in\whfK$.
That $\whAmax\subseteq\whL^*_0$ densely follows from
\[
\Amax\subseteq\Lmax=
P(L)^{1/2}\whL^*_0P(L)^{-1/2}\,.
\]
Thus $\Pi_{\whGm^\prime}$ is an isometric
boundary triple for $\whL^*_0$. Since
$u+\whk$ is the sum of a $\fH_2$-function
$u+\whk-\whg_z(c)$, $c=c(\whk)$, and
an eigenvector $\whg_z(c)$ of $\whL^*_0$,
it follows that $\whGm^\prime_0\subseteq\whGm_0$
and $\whGm^\prime_1\subseteq\whGm_1$.
Since moreover $\fN\dsum\whfK$ is a core for $\whL^*_0$,
this shows $\clo\whGm^\prime=\whGm$.

The last statement of the theorem uses
$\fN_z(\whAmax)=\fN_z(\whL^*_0)$.
\end{proof}
\begin{rems}
\begin{SL}
\item[1.]
$\cG_Z$ need not be Hermitian.
\item[2.]
Since $T_0$ is only essentially self-adjoint,
$\Pi_{\whGm^\prime}$ is not a
$B$-generalized boundary triple for $\whL^*_0$
\cite[Definition~1.5]{Derkach17}.
(In \cite[Lemma~5.5(ii)]{Derkach06} $A_0$ must be
self-adjoint.)
\end{SL}
\end{rems}

\end{appendix}
% \bibliography{peak-2021}

\begin{bibdiv}
\begin{biblist}

\bib{Albeverio10}{article}{
      author={Albeverio, S.},
      author={Cognola, G.},
      author={Spreafico, M.},
      author={Zerbini, S.},
       title={Singular perturbations with boundary conditions and the {C}asimir
  effect in the half space},
        date={2010},
     journal={J. Math. Phys.},
      volume={51},
       pages={063502},
}

\bib{Albeverio13}{article}{
      author={Albeverio, S.},
      author={Fassari, S.},
      author={Rinaldi, F.},
       title={A remarkable spectral feature of the {S}chr\"{o}dinger
  {H}amiltonian of the harmonic oscillator perturbed by an attractive
  $\delta^\prime$-interaction centred at the origin: double degeneracy and
  level crossing},
        date={2013},
     journal={J. Phys. A: Math. Theor.},
      volume={46},
      number={38},
       pages={385305},
}

\bib{Albeverio97}{article}{
      author={Albeverio, S.},
      author={Kurasov, P.},
       title={Rank {O}ne {P}erturbations of {N}ot {S}emibounded {O}perators},
        date={1997},
     journal={Integr. Equ. Oper. Theory},
      volume={27},
       pages={379\ndash 400},
}

\bib{Albeverio00}{book}{
      author={Albeverio, S.},
      author={Kurasov, P.},
       title={Singular {P}erturbations of {D}ifferential {O}perators},
      series={London Mathematical Society Lecture Note Series 271},
   publisher={Cambridge University Press, UK},
        date={2000},
}

\bib{Albeverio07}{article}{
      author={Albeverio, S.},
      author={Kuzhel, S.},
      author={Nizhnik, L.},
       title={Singularly perturbed self-adjoint operators in scales of
  {H}ilbert spaces},
        date={2007},
     journal={Ukrainian J. Math.},
      volume={59},
      number={6},
       pages={787\ndash 810},
}

\bib{Arlinski15}{incollection}{
      author={Arlinskii, Y.},
      author={Hassi, S.},
       title={${Q}$-functions and boundary triplets of nonnegative operators},
        date={2015},
   booktitle={Recent {A}dvances in {I}nverse {S}cattering, {S}chur {A}nalysis
  and {S}tochastic {P}rocesses},
      editor={Alpay, D.},
      editor={Kirstein, B.},
      series={Operator Theory: Advances and Applications},
      volume={244},
   publisher={Birkhauser},
       pages={89\ndash 130},
}

\bib{Azizov89}{book}{
      author={Azizov, T.},
      author={Iokhvidov, I.},
       title={Linear {O}perators in {S}paces with an {I}ndefinite {M}etric},
   publisher={John Wiley \& Sons. Inc.},
        date={1989},
}

\bib{Behrndt11}{article}{
      author={Behrndt, J.},
      author={Derkach, V.~A.},
      author={Hassi, S.},
      author={de~Snoo, H.},
       title={A realization theorem for generalized {N}evanlinna families},
        date={2011},
     journal={Operators and Matrices},
      volume={5},
      number={4},
       pages={679\ndash 706},
}

\bib{Behrndt20}{book}{
      author={Behrndt, J.},
      author={Hassi, S.},
      author={de~Snoo, H.},
       title={Boundary {V}alue {P}roblems, {W}eyl {F}unctions, and
  {D}ifferential {O}perators.},
   publisher={Birkhauser},
        date={2020},
}

\bib{Behrndt13}{article}{
      author={Behrndt, J.},
      author={Hassi, S.},
      author={de~Snoo, H.},
      author={Wietsma, R.},
      author={Winkler, H.},
       title={Linear fractional transformations of {N}evanlinna functions
  associated with a nonnegative operator},
        date={2013},
     journal={Compl. Anal. Oper. Theory},
      volume={7},
      number={2},
       pages={331\ndash 362},
}

\bib{Behrndt07}{article}{
      author={Behrndt, Jussi},
      author={Langer, Matthias},
       title={Boundary value problems for elliptic partial differential
  operators on bounded domains},
        date={2007},
     journal={J. Func. Anal.},
      volume={243},
       pages={536\ndash 565},
}

\bib{Behrndt18}{article}{
      author={Behrndt, Jussi},
      author={Langer, Matthias},
      author={Lotoreichik, Vladimir},
      author={Rohleder, Jonathan},
       title={Spectral enclosures for non-self-adjoint extensions of symmetric
  operators},
        date={2018},
     journal={J. Func. Anal.},
      volume={275},
      number={7},
       pages={1808\ndash 1888},
}

\bib{Bruning08}{article}{
      author={Br\"{u}ning, Jochen},
      author={Geyler, Vladimir},
      author={Pankrashkin, Konstantin},
       title={Spectra of self-adjoint extensions and applications to solvable
  {S}chr\"{o}dinger operators},
        date={2008},
     journal={Rev. Math. Phys.},
      volume={20},
      number={1},
       pages={1\ndash 70},
}

\bib{Derkach99}{article}{
      author={Derkach, V.},
       title={On generalized resolvents of {H}ermitian relations in {K}rein
  spaces},
        date={1999},
     journal={J. Math. Sci.},
      volume={97},
      number={5},
       pages={4420\ndash 4460},
}

\bib{Derkach03}{article}{
      author={Derkach, V.},
      author={Hassi, S.},
      author={de~Snoo, H.},
       title={Singular perturbations of self-adjoint operators},
        date={2003},
     journal={Math. Phys. Anal. Geom.},
      volume={6},
      number={4},
       pages={349\ndash 384},
}

\bib{Derkach06}{article}{
      author={Derkach, V.},
      author={Hassi, S.},
      author={Malamud, M.},
      author={de~Snoo, H.},
       title={Boundary relations and their {W}eyl families},
        date={2006},
     journal={Trans. Amer. Math. Soc.},
      volume={358},
      number={12},
       pages={5351\ndash 5400},
}

\bib{Derkach09}{article}{
      author={Derkach, V.},
      author={Hassi, S.},
      author={Malamud, M.},
      author={de~Snoo, H.},
       title={Boundary relations and generalized resolvents of symmetric
  operators},
        date={2009},
     journal={Russ. J. Math. Phys.},
      volume={16},
      number={1},
       pages={17\ndash 60},
}

\bib{Derkach91}{article}{
      author={Derkach, V.~A.},
      author={Malamud, M.~M.},
       title={Generalized {R}esolvents and the {B}oundary {V}alue {P}roblems
  for {H}ermitian {O}perators with {G}aps},
        date={1991},
     journal={J. Func. Anal.},
      volume={95},
      number={1},
       pages={1\ndash 95},
}

\bib{Derkach12}{incollection}{
      author={Derkach, Vladimir},
      author={Hassi, Seppo},
      author={Malamud, Mark},
      author={de~Snoo, Henk},
       title={Boundary triplets and {W}eyl functions. {R}ecent developments},
        date={2012},
   booktitle={Operator {M}ethods for {B}oundary {V}alue {P}roblems, {L}ondon
  {M}ath. {S}oc. {L}ecture {N}ote {S}eries},
      editor={Hassi, Seppo},
      editor={de~Snoo, Hendrik S.~V.},
      editor={Szafraniec, Franciszek~Hugon},
      volume={404},
   publisher={Cambridge University Press, UK},
       pages={161\ndash 220},
}

\bib{Derkach17}{article}{
      author={Derkach, Vladimir},
      author={Hassi, Seppo},
      author={Malamud, Mark~M.},
       title={Generalized boundary triples, {I}. {S}ome classes of isometric
  and unitary boundary pairs and realization problems for subclasses of
  {N}evanlinna functions},
        date={2020},
     journal={Math. Nachr.},
      volume={293},
      number={7},
       pages={1278\ndash 1327},
         url={https://onlinelibrary.wiley.com/doi/abs/10.1002/mana.201800300},
}

\bib{Dijksma05}{article}{
      author={Dijksma, A.},
      author={Kurasov, P.},
      author={Shondin, Yu.},
       title={High {O}rder {S}ingular {R}ank {O}ne {P}erturbations of a
  {P}ositive {O}perator},
        date={2005},
     journal={Integr. Equ. Oper. Theory},
      volume={53},
       pages={209\ndash 245},
}

\bib{Dijksma18}{article}{
      author={Dijksma, A.},
      author={Langer, H.},
       title={Compressions of self-adjoint extensions of a symmetric operator
  and {M}.{G}. {K}rein's resolvent formula},
        date={2018},
     journal={Integr. Equ. Oper. Theory},
      volume={90},
      number={41},
       pages={1\ndash 30},
}

\bib{Dijksma04a}{incollection}{
      author={Dijksma, A.},
      author={Langer, H.},
      author={Luger, A.},
      author={Shondin, Yu.},
       title={Minimal realizations of scalar generalized {N}evanlinna functions
  related to their basic factorization},
        date={2004},
   booktitle={Spectral {M}ethods for {O}perators of {M}athematical {P}hysics.
  {O}perator {T}heory: {A}dvances and {A}pplications},
      editor={Janas, J.},
      editor={Kurasov, P.},
      editor={Naboko, S.},
      volume={154},
   publisher={Birkh{\"a}user Basel},
}

\bib{Gallone20}{article}{
      author={Gallone, M.},
      author={Michelangeli, A.},
       title={Self-adjoint extensions with {F}riedrichs lower bound},
        date={2020},
     journal={Compl. Anal. Oper. Theory},
      volume={14},
      number={73},
       pages={23p},
}

\bib{Gohberg05}{book}{
      author={Gohberg, I.},
      author={Lancaster, P.},
      author={Rodman, L.},
       title={Indefinite {L}inear {A}lgebra and {A}pplications},
   publisher={Birkhauser},
        date={2005},
}

\bib{Halmos71}{article}{
      author={Halmos, P.~R.},
       title={Eigenvectors and adjoints},
        date={1971},
     journal={Lin. Alg. Appl.},
      volume={4},
       pages={11\ndash 15},
}

\bib{Hassi16}{article}{
      author={Hassi, S.},
      author={Wietsma, H.},
       title={Minimal realizations of generalized {N}evanlinna functions},
        date={2016},
     journal={Opuscula Math.},
      volume={36},
      number={6},
       pages={749\ndash 768},
}

\bib{Hassi09-b}{article}{
      author={Hassi, Seppo},
      author={Kuzhel, Sergii},
       title={On symmetries in the theory of finite rank singular
  perturbations},
        date={2009},
     journal={J. Func. Anal.},
      volume={256},
       pages={777\ndash 809},
}

\bib{Jursenas18a}{article}{
      author={Jur\v{s}\.{e}nas, R.},
       title={Computation of the unitary group for the {R}ashba spin--orbit
  coupled operator, with application to point-interactions},
        date={2018},
     journal={J. Phys. A: Math. Theor.},
      volume={51},
      number={1},
       pages={015203},
}

\bib{Jursenas21}{article}{
      author={Jur\v{s}\.{e}nas, R.},
       title={The {A}-model with mutually equal model parameters can lead to a
  {H}ilbert space model},
        date={2021},
     journal={Operators and Matrices},
      volume={15},
      number={4},
       pages={1319\ndash 1336},
}

\bib{Jursenas19}{article}{
      author={Jur\v{s}\.{e}nas, Rytis},
       title={On some extensions of the {A}-model},
        date={2020},
     journal={Opuscula Mathematica},
      volume={40},
      number={5},
       pages={569\ndash 597},
}

\bib{Kurasov03}{article}{
      author={Kurasov, P.},
       title={$\mathcal{H}_{-n}$-perturbations of self-adjoint operators and
  {K}rein's resolvent formula},
        date={2003},
     journal={Integr. Equ. Oper. Theory},
      volume={45},
      number={4},
       pages={437\ndash 460},
}

\bib{Kurasov03b}{article}{
      author={Kurasov, P.},
      author={Pavlov, Yu.~V.},
       title={On field theory methods in singular perturbation theory},
        date={2003},
     journal={Lett. Math. Phys.},
      volume={64},
      number={2},
       pages={171\ndash 184},
}

\bib{Kurasov09}{article}{
      author={Kurasov, Pavel},
       title={Triplet extensions {I}: {S}emibounded operators in the scale of
  {H}ilbert spaces},
        date={2009},
     journal={Journal d'Analyse Mathematique},
      volume={107},
      number={1},
       pages={252\ndash 286},
}

\bib{Langer77}{article}{
      author={Langer, H.},
      author={Textorius, B.},
       title={On generalized resolvents and {Q}-functions of symmetric linear
  relations (subspaces) in {H}ilbert space},
        date={1977},
     journal={Pacific. J. Math.},
      volume={72},
      number={1},
       pages={135\ndash 165},
         url={https://projecteuclid.org/euclid.pjm/1102811276},
}

\bib{Simon05a}{book}{
      author={Simon, B.},
       title={Trace {I}deals and {T}heir {A}pplications},
     edition={2},
      series={Mathematical Surveys and Monographs},
   publisher={American Mathematical Society},
        date={2005},
      volume={120},
}

\end{biblist}
\end{bibdiv}


\end{document}